\documentclass[11pt,psamsfonts]{amsart}
\usepackage{amsmath}
\usepackage{amsthm}
\usepackage{amssymb}
\usepackage{amscd}
\usepackage{amsfonts}
\usepackage{wrapfig}
\usepackage{amsbsy}
\usepackage{epsfig,afterpage}
\usepackage[dvips]{psfrag}
\usepackage[all]{xy}
\usepackage{pstcol}
\newcommand{\R}{\ensuremath{\mathbb{R}}}

\newcommand{\N}{\ensuremath{\mathbb{N}}}

\newtheorem {lemma}  {Lemma}
\newtheorem {definition}  {Definition}
\newtheorem {remark}  {Remark}
\newtheorem {example} {Example}

\newcommand{\cqd}{\begin{flushright} \vspace{-.6cm}$\Box$
\end{flushright}}

\definecolor{red}{rgb}{1.,0.,0.}
\definecolor{blue}{rgb}{0.,0.,1.}
\definecolor{pink}{rgb}{1.,0.75,0.8}

\begin{document}

\title[On three-parameter families of Filippov systems $-$ The
Fold-Saddle singularity] {On three-parameter families of Filippov systems $-$ The
Fold-Saddle singularity.}

\author[C.A. Buzzi, T. de Carvalho and M.A. Teixeira]
{Claudio A. Buzzi$^1$, Tiago de Carvalho$^2$ and\\ Marco A.
Teixeira$^3$}

\address{$^1$ IBILCE--UNESP, CEP 15054--000,
S. J. Rio Preto, S\~ao Paulo, Brazil}

\address{$^2$ FC--UNESP, 
Bauru, S\~ao Paulo, Brazil}

\address{$^3$ IMECC--UNICAMP, CEP 13081--970, Campinas,
S\~ao Paulo, Brazil}

\email{buzzi@ibilce.unesp.br}

\email{tcarvalho@fc.unesp.br}

\email{teixeira@ime.unicamp.br}


\date{}
\dedicatory{} \maketitle


\begin{abstract}

This paper presents  results concerning bifurcations of $2D$
piecewise-smooth vector fields. In particular, the generic
unfoldings of codimension three fold-saddle singularities of
Filippov systems, where a boundary-saddle and a fold coincide, are
considered and the bifurcation diagrams exhibited.

\end{abstract}

\section{Introduction}\label{secao introducao}

The general purpose of this paper is to study 
non-smooth vector fields (NSVF's for short), also called Filippov
systems, $Z$ 
represented by the following three-parameter
family of differential equations in $\R^2$:

\begin{equation}\label{eq fold-sela com parametros}
\overline{Z}^{\tau}_{\lambda , \mu, \beta} = \left\{
      \begin{array}{ll}
        X_{\lambda} =  \left(
              \begin{array}{c}
                  1 \\
                \alpha_1(\tau)(x - \lambda) + \alpha_2(\tau)(x-\lambda)^2
\end{array}
      \right)
 & \hbox{if $y \geq 0$,} \\
         Y_{\mu,\beta} = \left(
              \begin{array}{c}
                \frac{\mu}{2} x + \frac{(\mu - 2)}{2}(y + \beta) \\
                \frac{(\mu - 2)}{2} x + \frac{\mu}{2}(y + \beta)
\end{array}
      \right)& \hbox{if $y \leq 0,$}
      \end{array}
    \right.
\end{equation}
where  $\tau= \{ inv, vis \}$, $\alpha_1(inv)=-1$,
$\alpha_1(vis)=1$, $\alpha_2(inv)=1$, $\alpha_2(vis)=0$ and
$(\lambda,\beta,\mu) \in (-1,1) \times (-\sqrt{3}/2,\sqrt{3}/2)
\times (- \varepsilon_0, 1)$ with $\varepsilon_0 > 0$ sufficiently
small.

In Figure \ref{fig fold-sela} (respectively, Figure \ref{fig
fold-visivel-sela}) we consider $\lambda=\mu=\beta=0$ and $\tau=inv$
(respectively, $\tau=vis$) in Equation (\ref{eq fold-sela com
parametros}).

\begin{figure}[!h]
\begin{minipage}[b]{0.493\linewidth}
\begin{center}
 \epsfxsize=3.5cm \epsfbox{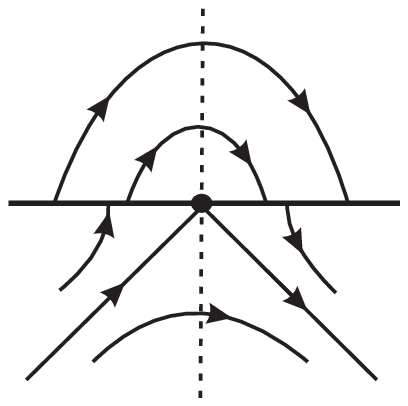}
\caption{\small{(Invisible) Fold-Saddle Singularity.}} \label{fig
fold-sela}
\end{center}
\end{minipage} \hfill
\begin{minipage}[b]{0.49\linewidth}
\begin{center}
 \epsfxsize=3.5cm \epsfbox{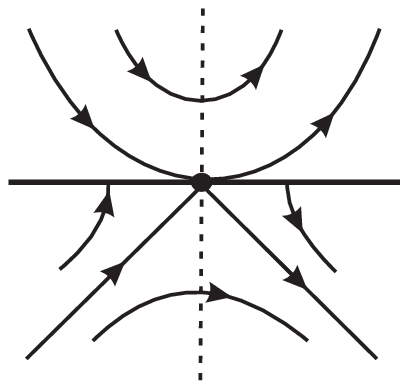}
\caption{\small{(Visible) Fold-Saddle Singularity.}} \label{fig
fold-visivel-sela}
\end{center}
\end{minipage}\end{figure}

Call $S_{X}= \{ X^2 = 0 \}$ and $S_Y=\{ Y^2=0 \}$ where
$X_{\lambda}=X=(X^1,X^2)$ and $Y_{\mu,\beta}=Y=(Y^1,Y^2)$. Denote
$\overline{Z}^{\tau}_{\lambda , \mu,
\beta}=(X_{\lambda},Y_{\mu,\beta})=(X,Y)$. In short our goal is to
study the local dynamics of  $\overline{Z}^{\tau}_{\lambda , \mu,
\beta}$ consisting of two smooth vector fields $X_{\lambda}$ and
$Y_{\mu,\beta}$ in $\R^2$ such that on one side of a smooth surface
$\Sigma = \{ y = 0 \}$ we take $\overline{Z}^{\tau}_{\lambda , \mu,
\beta} = X_{\lambda}$ and on the other side
$\overline{Z}^{\tau}_{\lambda , \mu, \beta} = Y_{\mu,\beta}$. We
mention two particular generic situations that occur in our system
when $\beta \neq 0$. The first one is the \textit{fold-fold
singularity} (or two-fold singularity $-$ see Definition
\ref{definicao fold-fold singularidade}), studied in some detail  in
\cite{Marcel} and  \cite{Kuznetsov} (among others), in which the
trajectories of both $X$ and $Y$ have a quadratic tangency to
$\Sigma$ at a point $p \in \Sigma$. The other situation is the
occurrence of branches of \textit{canard cycles} (see Definition
\ref{definicao canard cycle})  which are typical minimal sets that
appear in NSVF's. 

It is worth saying that the set $\Omega^r$ of all NSVF's $Z=(X,Y)$
as below described (with some specified topology) is a
differentiable manifold modeled in a Banach space (see Section
\ref{secao teoria basica} for details).

We now give a brief and rough overall description of the main
features of this work:
\begin{itemize}
\item In $\Omega^r$ all NSVF's presenting a generic fold-saddle singularity form a
codimension-three submanifold $W_1$ in such a way that any $Z \in
W_1$ is locally equivalent (see Definition \ref{definicao
sigma-equivalencia}) to the above described
$\overline{Z}^{\tau}_{\lambda , \mu, \beta}$ with
$\lambda=\mu=\beta=0$.


\item  The bifurcation diagram of $\overline{Z}^{\tau}_{\lambda , \mu, \beta}$ is exhibited for $\tau=inv$ and $\tau =vis$.

\end{itemize}

We emphasize that this paper is inserted in a larger 2D
classification program where are included papers like \cite{T1, T,
Kuznetsov, Marcel}. One of the goals of this program is to classify
(via topological equivalence) typical singularities of NSVF's. For
this purpose it is necessary to present generic unfoldings and give
non-degeneracy conditions on the system in order to characterize the
codimension of the singularity. The bifurcation diagram of the
codimension-three singularity presented here includes that one
exhibited in \cite{Marcel}. The later claim is commented in the
sequel. We finish this introduction presenting a physical model
where a fold-saddle singularity can be found.

%

\begin{example}\label{exemplo}One of the reasons that second order linear equations with
constant coefficients are worth studying is that they serve as
mathematical models of some important physical processes. Two
important areas of application are in the fields of mechanical and
electrical oscillations. For example, the motion of a mass on a
vibrating spring, the torsional oscillations of a shaft with a
flywheel, the flow of electric current in a simple series circuit,
and many other physical problems are all described by the solution
of an equation of the form
\begin{equation}\label{eq1}
a \ddot{x} + b \dot{x} + c x = g ( t,x ).
\end{equation}
Here we consider the external force $g$ not depending on $t$ but
depending on $x$. For example consider
\begin{equation}\label{eq2}
g(x)=Ax+1-sgn(x) \mbox{ with }A>\frac{c}{a}.
\end{equation}
Now if we call $\dot{x}=y$ then the equation \eqref{eq1} with $g$
given by \eqref{eq2} became
\begin{equation}\label{eq3} \dot{x}=y, \quad \dot{y}=\left(A-\frac{c}{a}\right)x -\frac{b}{a}y,\quad \mbox{if }x>0\end{equation} and
\begin{equation}\label{eq4}\dot{x}=y, \quad \dot{y}=\left(A-\frac{c}{a}\right)x -\frac{b}{a}y+2,\quad \mbox{if }x<0.\end{equation}
System \eqref{eq3} has a saddle equilibrium at the origin because
$A>\displaystyle\frac{c}{a}$. And system \eqref{eq4} has an
invisible fold (see Definition \ref{definicao fold point}) at the
origin. If $g(x)=Ax-1+sgn(x) \mbox{ with
}A>\displaystyle\frac{c}{a}$ then it has a saddle equilibrium and a
visible fold.\cqd
\end{example}

\vspace{.3cm}

 The paper is organized as follows: In Section \ref{secao
enunciado dos teoremas} we present the main results of the paper. In
Section \ref{secao teoria basica} we give some basic concepts about
NSVF's in order to setting the problem in Section \ref{secao setting
the problem}.  The remaining sections are dedicated to prove the
main results of the paper.

\section{Statement of the Main Results}\label{secao enunciado dos
teoremas}
%
In what follows consider \begin{equation}\label{eq mu zero}
\mu_0(\beta) = 2 - (12 \beta/(-3 + 6 \beta + \sqrt{9- 12 \beta^2})).
\end{equation}\vspace{.5cm}


\noindent {\bf Theorem 1.} {\it Assume $\tau=inv$ and
$\mu=\mu_0(\beta)$ in Equation (\ref{eq fold-sela com parametros}). 
 Then its bifurcation diagram in the
$(\lambda,\beta)-$plane contains $19$ distinct phase portraits (see
Figure \ref{fig diagrama bif teo 1}).}

\vspace{.5cm}


First of all observe that if $\beta >0$ and $\lambda = -1/2 +
\sqrt{9 - 12 \beta^2}/6$ in Theorem 1, then the NSVF has a
homoclinic loop surrounding a non hyperbolic singularity. So, it is
easy to see that the cases covered by Theorem 1 do not represent the
full unfolding of the (invisible) fold-saddle singularity and the
next two theorems become necessary.


\vspace{.5cm}

\noindent {\bf Theorem 2.} {\it Assume $\tau=inv$ and $\mu_0(\beta)
<\mu < 1$ in Equation (\ref{eq fold-sela com parametros}). Then its
bifurcation diagram in the $(\lambda,\beta)-$plane contains
 $21$ distinct phase portraits
 (see Figure \ref{fig diagrama bif teo 2}).}

\vspace{.5cm}

\noindent {\bf Theorem 3.} {\it Assume $\tau=inv$ and $-
\varepsilon_0<\mu < \mu_0(\beta)$ in Equation (\ref{eq fold-sela com
parametros}). Then its bifurcation diagram in the
$(\lambda,\beta)-$plane contains $21$  distinct phase portraits (see
Figure \ref{fig diagrama bif teo 2}).}

\vspace{.5cm}

\begin{figure}[!h]
\begin{center}\psfrag{a}{$(a)$} \psfrag{b}{$(b)$} \psfrag{c}{$(c)$}
\psfrag{0}{$13_3$} \psfrag{2}{$14_3$} \psfrag{0}{$12_3$}
\epsfxsize=12.5cm \epsfbox{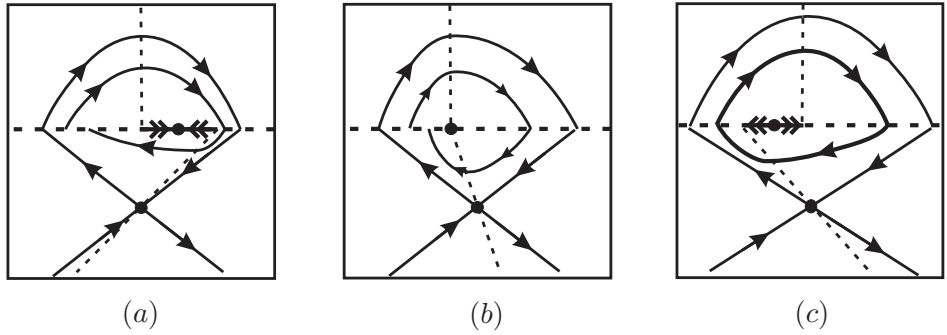} \caption{\small{A like Hopf
bifurcation. The singularities in (a) and (c) (that appear in
Theorems 2 and 3 respectively) are called $\Sigma-$attractor and
$\Sigma-$repeller throughout this paper (see Figure \ref{fig
singularidades} cases (b) and (c)). In (b) we observe a non
hyperbolic singularity present in Theorem 1.}} \label{fig hopf}
\end{center}
\end{figure}

\begin{remark}
The bifurcation diagrams exhibited in Theorems 2 and 3 present a
homoclinic loop surrounding a hyperbolic singularity. Observe that,
under the conditions of Theorem 2 (respectively, Theorem 3) this
singularity is an attractor (respectively, repeller) as illustrated
in Figure \ref{fig hopf}, cases $(a)$ and $(b)$. Moreover, when the
parameter $\mu$ varies from $- \varepsilon_0$ to $1$ there is an
element $\mu_{0}(\beta) \in (- \varepsilon_0,1)$, given by Equation
(\ref{eq mu zero}), such that $\overline{Z}^{inv}_{\lambda ,
\mu_{0}(\beta), \beta}$ presents a \textit{like Hopf bifurcation
phenomenon} as illustrated in Figure \ref{fig hopf}. This phenomenon
is fully treated in \cite{Marcel} and \cite{Kuznetsov} and it is not
covered by Theorems 1, 2 and 3.
\end{remark}

\begin{remark}
In Theorems 2 and 3 beyond the cases $(a)$ and $(c)$ in Figure
\ref{fig hopf}, it also appears another relevant phenomenon that is
the presence of a fold-fold singularity. Note that in Theorem 1 it
occurs simultaneously the occurrence of a loop (a global phenomenon)
and a fold-fold singularity (a local phenomenon).
\end{remark}

\noindent {\bf Theorem 4.} {\it Assume $\tau=vis$ in Equation
(\ref{eq fold-sela sem variar}) or equivalently, take $\tau=vis$ and
$\mu=0$ in Equation (\ref{eq fold-sela com parametros}). Then its
bifurcation diagram in the $(\lambda,\beta)-$plane contains $13$
distinct phase portraits (see Figure \ref{fig diagrama bif teo 4 5 e
6}).}

\vspace{.5cm}

The cases covered by Theorem 4 do not represent the full unfolding
of the (visible) fold-saddle singularity. So, the next two theorems
are necessary. Each one of them describes a distinct generic
codimension two singularity.\\

 \noindent {\bf Theorem 5.} {\it Assume $\tau=vis$ and $0 <\mu < 1$ in
Equation (\ref{eq fold-sela com parametros}). Then its bifurcation
diagram in the $(\lambda,\beta)-$plane contains $13$ distinct phase
portraits
 (see Figure \ref{fig diagrama bif teo 4 5 e
6}).}

\vspace{.5cm}

 \noindent {\bf Theorem 6.} {\it Assume $\tau=vis$ and $\varepsilon_0 <\mu < 0$
in Equation (\ref{eq fold-sela com parametros}). Then its
bifurcation diagram in the $(\lambda,\beta)-$plane contains $13$
distinct phase portraits
 (see Figure \ref{fig diagrama bif teo 4 5 e
6}).}

\begin{remark}
In Theorems 5 and 6 one observes the birth of a singularity on
$\Sigma$. This singularity behaves like a saddle (see Figure
\ref{fig singularidades}$-(a)$) and is known as $\Sigma-$saddle.
This phenomenon does not occur under the conditions presented in
Theorem 4.
\end{remark}

\section{Basic Theory about NSVF}\label{secao teoria basica}

Let $K\subseteq \R ^{2}$ be a compact set  such that $\partial K$ is
a smooth curve. Consider $\Sigma \subseteq K$ given by $\Sigma
=f^{-1}(0),$ where $f:K\rightarrow \R$ is a smooth function having
$0\in \R$ as a regular value (i.e. $\nabla f(p)\neq 0$, for any
$p\in f^{-1}({0}))$ such that $\partial K \cap \Sigma = \emptyset$
or $\partial K \pitchfork \Sigma$. Clearly the \textit{switching
manifold} $\Sigma$ is the separating boundary of the regions
$\Sigma_+=\{q\in K \, | \, f(q) \geq 0\}$ and $\Sigma_-=\{q \in K \,
| \, f(q)\leq 0\}$. We can assume that $\Sigma$ is represented,
locally around a point $q=(x,y)$, by the function $f(x,y)=y.$\\

Designate by $\chi^r$ the space of $C^r-$vector fields on $K$
endowed with the $C^r-$topology with $r\geq 1$ large enough for our
purposes. Call \textbf{$\Omega^r=\Omega^r(K,f)$} the space of vector
fields $Z: K  \rightarrow \R ^{2}$ such that
$$
 Z(x,y)=\left\{\begin{array}{l} X(x,y),\quad $for$ \quad (x,y) \in
\Sigma_+,\\ Y(x,y),\quad $for$ \quad (x,y) \in \Sigma_-,
\end{array}\right.
$$
where $X=(X^1,X^2)$ and $Y = (Y^1,Y^2)$ are in $\chi^r.$ We write
$Z=(X,Y),$ which we will accept to be multivalued in points of
$\Sigma.$ The trajectories of $Z$ are solutions of  ${\dot q}=Z(q),$
which has, in general, discontinuous righthand side. The basic
results of differential equations, in this context, were stated by
Filippov in \cite{Fi}. Related theories can be found in \cite{K, T} among others.\\

\subsection{Orbits, trajectories and singularities of NSVF's}\label{subsecao orbitas trajet e
sing}

In what follows we will use the notation
\[X.f(p)=\left\langle \nabla f(p), X(p)\right\rangle \quad \mbox{
and } \quad X^i.f(p)=\left\langle \nabla (X^{i-1}. f)(p),
X(p)\right\rangle, i\geq 2
\]
where $\langle . , . \rangle$ is the usual inner product in
$\R^2$.\\

Following the Filippov rule, we  distinguish the following regions
on the discontinuity set $\Sigma:$

\begin{itemize}
\item   \textbf{Crossing region:}   $\Sigma_c = \{ p \in \Sigma \, |
\, (X.f(p))(Y.f(p))>0 \}$.

\item    \textbf{Escaping region:}   $\Sigma_e =  \{ p \in \Sigma \, |
\, X.f(p) > 0 \mbox{ and } Y.f(p) < 0 \}$.

\item    \textbf{Sliding region:}   $\Sigma_s =  \{ p \in \Sigma \, |
\, X.f(p) < 0 \mbox{ and } Y.f(p) > 0 \}$.
\end{itemize}

\begin{figure}[h]
\begin{center}
\psfrag{A}{$q$} \psfrag{B}{$q + Y(q)$} \psfrag{C}{$q + X(q)$}
\psfrag{D}{} \psfrag{E}{\hspace{1cm}$Z^\Sigma(q)$}
\psfrag{F}{\hspace{-.35cm}$\Sigma_s$} \psfrag{G}{} \epsfxsize=5.5cm
\epsfbox{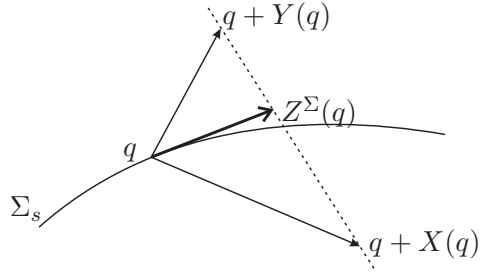} \caption{\small{Filippov's
convention.}} \label{fig def filipov}
\end{center}
\end{figure}

Consider $Z=(X,Y) \in \Omega^r$ and $p \in \Sigma_e \cup \Sigma_s$.
The \textit{sliding vector field} $Z^{\Sigma}$  associated to $Z$ at
$p$ is the convex combination of $X(p)$ and $Y(p)$ tangent to
$\Sigma$ at $p$ (see Figure \ref{fig def filipov}).


We say that $q\in\Sigma$ is a \textit{$\Sigma-$regular point} if
\begin{itemize}
\item $(X.f(q))(Y.f(q))>0$
or
\item $(X.f(q))(Y.f(q))<0$ and $Z^{\Sigma}(q)\neq0$ (that is $q\in\Sigma_e\cup\Sigma_s$ and it is not an equilibrium
point of $Z^{\Sigma}$).\end{itemize}

The points of $\Sigma$ which are not $\Sigma-$regular are called
\textit{$\Sigma-$singular}. We distinguish two subsets in the set of
$\Sigma-$singular points: $\Sigma^t$ and $\Sigma^p$. Any $q \in
\Sigma^p$ is called a \textit{pseudo equilibrium of $Z$} and it is
characterized by $Z^{\Sigma}(q)=0$. Any $q \in \Sigma^t$ is called a
\textit{tangential singularity} and is characterized by
$Z^{\Sigma}(q) \neq 0$ and
$(X.f(q))(Y.f(q)) =0$.\\

%

\begin{figure}[h]
\begin{center}
\psfrag{a}{(a)} \psfrag{b}{(b)} \psfrag{c}{(c)} \psfrag{D}{}
\psfrag{E}{\hspace{1cm}$Z^\Sigma(q)$}
\psfrag{F}{\hspace{-.35cm}$\Sigma_s$} \psfrag{G}{} \epsfxsize=11cm
\epsfbox{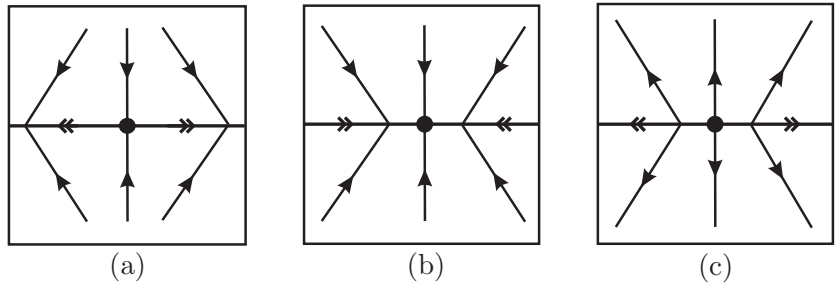} \caption{\small{Pseudo equilibria:
(a) $\Sigma-$saddle, (b) $\Sigma-$attractor and (c)
$\Sigma-$repeller.}} \label{fig singularidades}
\end{center}
\end{figure}

A pseudo equilibrium $q \in \Sigma^p$ is a \textit{$\Sigma-$saddle}
provided that one of the following condition is satisfied: (i)
$q\in\Sigma_e$ and $q$ is an attractor for $Z^{\Sigma}$ or (ii)
$q\in\Sigma_s$ and $q$ is a repeller for $Z^{\Sigma}$. A pseudo
equilibrium $q\in\Sigma^p$ is a $\Sigma-$\textit{repeller} (resp.
$\Sigma-$\textit{attractor}) provided $q\in\Sigma_e$ (resp. $q \in
\Sigma_s$) and $q$ is a repeller (resp. attractor) equilibrium point
for
$Z^{\Sigma}$ (see Figure \ref{fig singularidades}). 

\begin{definition}\label{definicao fold point} We say that
$p_0 \in \Sigma^t$ is a \textbf{$\mathbf{\Sigma-}$fold point} of $X
\in \chi^r$ if $X.f(p_0)=0$ but $X^{2}.f(p_0)\neq0.$ Moreover,
$p_0\in\Sigma$ is a \textit{visible} (respectively {\it invisible})
$\Sigma-$fold point of $X$ if $X.f(p_0)=0$ and $X^{2}.f(p_0)> 0$
(respectively $X^{2}.f(p_0)< 0$).\end{definition}

\begin{definition}\label{definicao fold-fold singularidade}
Let $Z=(X,Y) \in \Omega^r$. We say that $p \in \Sigma^t$ is a
\textbf{fold-fold singularity} of $Z$ if $p$ is a $\Sigma-$fold
point for both $X$ and $Y$.
\end{definition}

\begin{definition}\label{definicao fold-sela singularidade}
Let $Z=(X,Y) \in \Omega^r$. We say that $q \in \Sigma$ is a
\textbf{fold-saddle singularity} of $Z$ if $q$ is a $\Sigma-$fold
point of $X$ and a saddle equilibrium of $Y$ (in this case $q$ is
called a boundary-saddle of $Y$).
\end{definition}

The following construction is presented  in
\cite{Marco-enciclopedia}. Let $Z=(X,Y) \in \Omega^r$ such that $T$
is an invisible $\Sigma-$fold point of $X$. From Implicit Function
Theorem, for each $p \in \Sigma$ in a neighborhood $\mathcal{V}_{T}$
of $T$ we derive that there exists a unique $t(p)$ such that the
orbit$-$solution $t \mapsto \varphi_X(t,p)$ of $X$ through $p$ meets
$\Sigma$ at a point $\widetilde{p}=\varphi_X(t(p),p)$. Define the
 map $\gamma_X:\mathcal{V}_{T} \cap \Sigma \rightarrow
\mathcal{V}_{T} \cap \Sigma$ by $\gamma_X(p)=\widetilde{p}$. This
map is a $C^{r}-$diffeomorphism and satisfies: $\gamma_{X}^{2}=Id$.
Analogously, when $\widetilde{T}$ is an invisible $\Sigma-$fold
point of $Y$  we define the  map
$\gamma_Y:\mathcal{V}_{\widetilde{T}} \cap \Sigma \rightarrow
\mathcal{V}_{\widetilde{T}} \cap \Sigma$  associated to $Y$ which
satisfies $\gamma_{Y}^{2}=Id$. We define now the first return map
associated to $Z=(X,Y)$:

\begin{definition}\label{aplicacao-primeiro-retorno}
The \textbf{first return map} $\varphi_Z: \mathcal{T} \rightarrow
\mathcal{T}$ is defined by the composition $\varphi_Z=\gamma_Y \circ
\gamma_X$ when both $\gamma_X$ and $\gamma_Y$ are well defined in
$\mathcal{T} \subset \Sigma$.\end{definition}

\begin{remark} $\varphi_Z$ is an area-preserving map.\\  \end{remark}

\begin{definition}\label{definicao canard cycle}
A curve $\Gamma$ is a \textbf{canard cycle} of $Z=(X,Y) \in
\Omega^r$ if it is closed and composed by orbit-arcs of at least two
of the vector fields $X |_{\Sigma_{+}}$, $Y |_{\Sigma_{-}}$ and
$Z^{\Sigma}$. We say that $\Gamma$ is \textbf{hyperbolic} if
$\varphi_Z '(p) \neq 1$, where $\varphi_Z$ is the first return map
defined on a segment $T$ with $p \in T \pitchfork \gamma$.
\end{definition}

\begin{definition}\label{definicao grafico} Consider $Z \in \Omega^r$. A closed path $\Delta$ is
a \textbf{$\mathbf{\Sigma-}$graph} if it is a union of equilibria,
pseudo equilibria, tangential singularities of $Z$ and orbit-arcs of
$Z$ joining these points in such a way that $\Delta \cap \Sigma \neq
\emptyset$.\\ 
\end{definition}

\begin{figure}[h!]
\begin{minipage}[b]{0.4\linewidth}
\begin{center}
\epsfxsize=4.7cm \epsfbox{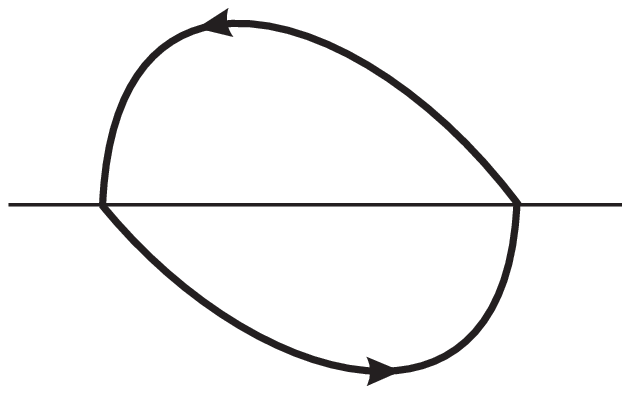} \caption{\small{Canard
cycle.}} \label{fig canard I}
\end{center}
\end{minipage} \hfill
\begin{minipage}[b]{0.5\linewidth}
\begin{center}\psfrag{4}{$12_1$} \epsfxsize=4.5cm
\epsfbox{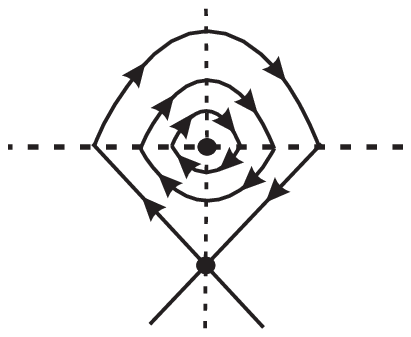} \caption{\small{$\Sigma-$graph
surrounding a $\Sigma-$center.}} \label{fig sigma centro}
\end{center}
\end{minipage}
\end{figure}

\begin{definition}\label{definicao centro}
Consider $Z \in \Omega^r$. A point $q \in \Sigma$ is a
\textbf{$\mathbf{\Sigma-}$center} if there is a neighborhood $U$ of
$q$ filled up with a one-parameter family $\gamma_{U}$ of canard cycles of $Z$ in such a way that $\gamma_{U} \cap \Sigma \subset \Sigma_c$ (see Figure \ref{fig sigma centro}). 
\end{definition}


\subsection{Structural Stability on $\Omega^r$}\label{subsecao estabilidade estrutural}

Bifurcation theory describes how continuous variations of parameter
va\-lues in a  dynamical system can, through topological changes,
cause the phase portrait to change suddenly. In this paper we focus
on  certain structurally unsta\-ble NSVF's within a generic context.
In \cite{Andronov} the concept of $k^{\mbox{\small th}}-$order
structural stability is presented; in a local approach such setting
gives rise to the notion of a codimension $k$ singularity. Now we
present the concept of equivalence which will guide us for all the
paper.


\begin{definition}\label{definicao sigma-equivalencia} Two NSVF's $Z, \, \widetilde{Z} \in
\Omega^r(K,f)$ defined in open sets $U, \, \widetilde{U} \subset K$
with switching manifolds $\Sigma\subset U$ and $\widetilde{\Sigma}
\subset \widetilde{U}$ respectively are
\textbf{$\mathbf{\Sigma-}$equivalent} if there exists an orientation
preserving ho\-meo\-mor\-phism $h:U \rightarrow \widetilde{U}$ which
sends $\Sigma$ to $\widetilde{\Sigma}$ and sends orbits of $Z$
(respectively $Z^{\Sigma}$) to orbits of $\widetilde{Z}$
(respectively $\widetilde{Z}^{\widetilde{\Sigma}}$). From this
definition the concept of \textbf{local structural stability} in
$\Omega^r$ is naturally obtained.
\end{definition}

As we said in Section \ref{secao introducao} our paper is a
generalization of some papers that unfold typical singularities of
NSVF's. Below we present the program used in the literature and in our paper to exhibit the
diagram bifurcation of a singularity of NSVF's.\\

Let $Z \in \Omega^r$ and $p \in \Sigma$. Following the approach in
\cite{Soto} (and also exposed in \cite{Marcel}, Subsection 3.1, pg
1982), we get that:

\begin{itemize}
\item By Theorem 3.5 of \cite{Marcel} we already know the characterization
of the set $\Phi_{0}^{p} = \{ Z \in \Omega^r \, | \, Z \mbox{ is
locally structurally stable at } p \}$. In fact, $\Phi_{0}^{p}$ is
open and dense in $\Omega^r$. So, $\Phi_{0}^{p}$ is the codimension
zero local bifurcation set.


\item Let $\Omega_1 = \Omega^r \backslash \Phi_{0}^{p}$ and $\Phi_{1}^{p}
= \{ Z \in \Omega_1 \, | \, Z \mbox{ is locally structurally stable}
 \mbox{ at } p \mbox{ relative to  } \Omega_1 \}$. The set
$\Phi_{1}^{p}$ is the codimension $1$ local bifurcation set in
$\Omega^r$. The characterization $\Phi_{1}^{p}$ was given in
\cite{Kuznetsov}.

In addition, if $Z_0 \in \Phi_{1}^{p}$, it is also known (see
\cite{Marco-bollettino})  that there exists a neighborhood
$\mathcal{U}$  of $Z_0$ in $\Omega^r$ such that:

%

\hspace{1cm}$\mathbf{\circ}$ There exists a $C^r-$function $L :
\mathcal{U} \rightarrow \R$, satisfying $L(Z_0)=0$ and $DL_{Z_{0}}$,
the differential of $L$ at $Z_0$, is surjective. Moreover,
$L^{-1}(0)= \Phi_{1}^{p} \cap \mathcal{U}$.

\hspace{1cm}$\circ$ Consider now all the embeddings $\Theta :  ( -
\varepsilon , \varepsilon) \times \mathcal{U}  \rightarrow
\mathcal{U}$ transversal to $\Phi_{1}^{p}$ at some $Z \in
\Phi_{1}^{p}$ and such that $\Theta(0, Z_0)=Z$. We refer to such
$\Theta$ as an \textit{unfolding} of
$Z_0$. 

\item We consider now the set $\Omega_2 = \Omega_1
\backslash \Phi_{1}^{p}$ and similar objects $\Phi_{2}^{p}$ (the set
of codimension $2$ singularities) and families of objects $L :
\mathcal{U} \rightarrow \R^2$, with surjective derivative at $Z_0$
and embeddings $\Theta : (-\varepsilon , \varepsilon) \times (-
\zeta, \zeta) \times \mathcal{U} \rightarrow \mathcal{U}$.

\item In this way we can get sequences
$\Omega_k$ and $\Phi_{k}^{p}$ in $\Omega^r$, that establish a
program to characterize all codimension $k$ singularities.

\end{itemize}

\begin{definition}\label{definicao desdobramento generico}
Let $\mathcal{V}(0,\R^k)$ and $\mathcal{S}(0,\R^l)$ be neighborhoods
of $0$ in $\R^k$ and $\R^l$ respectively and let $\mathcal{U}$ be a
neighborhood of $Z_0$ in $\Omega^r$. We say that two
\textbf{unfoldings} $\Theta:\mathcal{V}(0,\R^k) \times \mathcal{U}
\rightarrow \mathcal{U}$ and $\Xi:\mathcal{S}(0,\R^l)
\times\mathcal{U} \rightarrow \mathcal{U}$ are \textbf{equivalent}
if there is a homomorphism  $A: \mathcal{V}(0,\R^k) \rightarrow
\mathcal{S}(0,\R^l)$ such that $A(\lambda) = \mu$ and for each $Z
\in \mathcal{U}$ the vector fields $\Theta(\lambda,Z)$ and
$\Xi(A(\lambda),Z)$ are $\Sigma-$equivalent according to
Definition \ref{definicao sigma-equivalencia}. 
Moreover, we say that an unfolding $\Theta(\lambda_{0},.)$ is a
\textbf{generic unfolding} if there is a neighborhood
$\mathcal{W}(\Theta(\lambda_{0},.))$ of $\Theta(\lambda_{0},.)$ such
that any unfolding $\Theta(\lambda,.) \in
\mathcal{W}(\Theta(\lambda_{0},.))$ is equivalent to
$\Theta(\lambda_{0},.)$.
\end{definition}

\begin{remark}\label{obs Marco desdobramento} In Definition \ref{definicao desdobramento
generico} it is important to say that the homomorphism $A$ does not
vary, necessarily, continuously with respect to the parameter
$\lambda \in \mathcal{V}(0,\R^k)$.
\end{remark}

\subsection{The Direction Function $H$}\label{subsecao funcao H}

Here we introduce a function that will be very usefull in the
sequel.\\

\begin{figure}[!h]
\begin{center}
\psfrag{A}{$A$} \psfrag{B}{$B$} \psfrag{C}{$C$} \psfrag{D}{$\Sigma$}
\psfrag{E}{$X$} \psfrag{F}{$Y$} \psfrag{G}{$C+Y(C)$} \psfrag{H}{$C+
X(C)$} \psfrag{I}{$p(C)$} \epsfxsize=5cm
\epsfbox{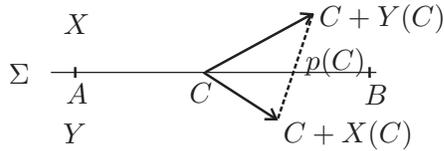} \caption{\small{Direction function.}}
\label{fig funcao direcao}
\end{center}
\end{figure}

In $(A,B)\subset \Sigma_e \cup \Sigma_s$, consider  the point
$C=(C_{1},C_2)$, the vectors $X(C)=(D_1,D_2)$ and $Y(C)=(E_1,E_2)$
(as illustrated in Figure \ref{fig funcao direcao}). The straight
segment passing through $C+X(C)$ and $C + Y(C)$ meets $\Sigma$ in a
point $p(C)$. We define the C$^r-$map
$$
\begin{array}{cccc}
  p: & (A,B) & \longrightarrow & \Sigma \\
     & z & \longmapsto & p(z).
\end{array}
$$
Since $\Sigma$ is the $x-$axis, we have that $C=(C_1,0)$ and $p(C)
\in \R \times \{ 0 \}$ can be identified with points in $\R$.
According with this identification, the \textit{direction function}
on $\Sigma$ is defined by
$$
\begin{array}{cccc}
  H: & (A,B) & \longrightarrow & \R \\
     & z & \longmapsto & p(z) - z.
\end{array}
$$

We obtain that $H$ is a C$^r-$map  and

\begin{itemize}
\item if $H(C) < 0$ then the orientation of $Z^{\Sigma}$ in a small neighborhood of $C$ is from $B$ to $A$;

\item if $H(C) = 0$ then  $C \in \Sigma^p$;

\item if $H(C) > 0$ then the orientation of $Z^{\Sigma}$ in a small neighborhood of $C$ is from $A$ to $B$.
\end{itemize}

Simple calculations show that $p(C_1) = \frac{E_2 (D_1+C_1)  - D_2
(E_1+C_1)}{E_2 - D_2}$ and consequently,
\begin{equation}\label{eq H}
H(C_1) = \frac{E_2 D_1  - D_2 E_1}{E_2 - D_2}.
\end{equation}

\begin{remark}\label{obs H positiva ou negativa nas dobras}
If $X.f(p)=0$ and $Y.f(p) \neq 0$ then, in a neighborhood $V_p$ of
$p$ in $\Sigma$, the direction function $H$ has the same sign of
$D_1$, where $X(p)=(D_1,D_2)$. In fact, $X.f(p)=0$ and $Y.f(p) \neq
0$ are equivalent to $D_2=0$ and $E_2 \neq 0$ in  (\ref{eq H}). So,
$\displaystyle\lim_{(D_2,E_2)\rightarrow(0,k_0)} H(p_1) = D_1$,
where $k_0 \neq 0$ and $p=(p_1,p_2)$.\\
\end{remark}


\section{Setting the problem}\label{secao setting the problem}

%
%
%
%
%
%
%
%
%
%
%
%
%
%
%

Let $\Gamma_{+} = \{ X|_{\Sigma_+} \mbox{ with } X \in \chi^r \}$
(respectively, $\Gamma_{-} = \{ X|_{\Sigma_-} \mbox{ with } X \in
\chi^r \}$). This means that $\Gamma_+$ (respectively, $\Gamma_-$)
is identified with $\chi^r$.

 Let
$\Gamma^{F}_{\Sigma_{+}} \subset \Gamma_+$ be the set of all
elements $X \in \Gamma_+$ having a $\Sigma-$fold point. 
 $\Gamma^{F}_{\Sigma_{+}}$ is an open set in $\Gamma_+$ (see
 \cite{T1}).
We may consider $f(x,y)=y$ and the following generic normal forms
$X_0(x,y)= (\alpha_1, \beta_1 x)$ with $\alpha_1=\pm1$ and $\beta_1=
\pm1$ (see \cite{vishik}, Theorem $2$).

Let $\Gamma^{S}_{\Sigma_{-}} \subset \Gamma_-$ be the set of all
elements $Y \in \Gamma_-$ presenting a hyperbolic saddle equilibrium
$S_{Y_0}$ on $\Sigma$ (called \textit{boundary saddle} in the
literature $-$ as a reference, in \cite{RR}, are exhibited some
border collisions in three-dimension NSVF's) and such that the
eigenspaces of
$DY_0(S_{Y_0})$ are transverse to $\Sigma$ at $S_{Y_0}$. 
$\Gamma^{S}_{\Sigma_{-}}$ is a codimension one submanifold of
$\Gamma_-$. Note that, $\Gamma^{S}_{\Sigma_{-}}$ has a $C^r
-$structure (for more details about this construction see
\cite{Dumortier-livro}). In fact, since $Y_0$ has a hyperbolic
boundary saddle $S_{Y_0}$ and the eigenspaces of $DY_0(S_{Y_0})$ are
transverse to $\Sigma$ at $S_{Y_0}$, then there exists a
neighborhood $\mathcal{V}(Y_0)$ of $Y_0$ in $\chi^r$ such that all
$Y \in \mathcal{V}(Y_0)$ have a hyperbolic boundary saddle
$S_Y=(x^*,y^*)$ with the same properties of the eigenspaces.
Moreover, the correspondence $L : \mathcal{V}(Y_0) \rightarrow
\R^2$, given by $L(Y)=S_Y$, is $C^r$ at $Y_0$. Define $\Pi :
\mathcal{V}(Y_0) \rightarrow \R$, given by $\Pi(Y)=(f \circ L)(Y)$.
We say that two vector fields $Y, \widetilde{Y} \in
\Gamma^{S}_{\Sigma_{-}}$ defined in open sets $U$ and
$\widetilde{U}$, respectively, are $C^{0}-$orbitally equivalent if
there exists an orientation preserving homeomorphism $h: U
\rightarrow \widetilde{U}$ that sends orbits of $Y$ to orbits of
$\widetilde{Y}$. From \cite{T1} we know that any $Y \in
\Gamma^{S}_{\Sigma_{-}}$ is generically $C^{0}-$orbitally
equi\-va\-lent to its linear part by a $\Sigma-$preserving
homeomorphism. And the linear saddle with eigenspaces transverse to
the $x-$axis has the generic normal forms $Y_0(x,y)= (\alpha_2 y,
\alpha_2 x)$ with $\alpha_2=\pm1$. So it is easy to see that the
generic unfolding of the singularity is given by $Y_{\beta}= (
\alpha_2 (y + \beta) , \alpha_2 x)$ where $\beta \in \R$.

%


Moreover: \begin{itemize}

\item There exists a $C^r -$function $\Pi: \mathcal{V}(Y_0)
\rightarrow \R$, such that $D \Pi_{Y_0}$ is surjective.

\item The correspondence $Y \mapsto S_Y$ is $C^r$, where $S_Y$ is a
saddle point of $Y$.

\item If $\Pi(Y) <0$ then $S_Y \in \Sigma_-$.

\item If $\Pi(Y) =0$ then $S_Y \in \Sigma$.

\item If $\Pi(Y) >0$ then $S_Y  \in \Sigma_+$.

\end{itemize}

In this paper we are concerned with the bifurcation diagram of
systems  $Z_0=(X_0,Y_0)$ in $\Omega^r$ such that $X_0 \in
\Gamma_{\Sigma_{+}}^{F}$, $Y_0 \in \Gamma_{\Sigma_{-}}^{S}$ and
$p_0=S_{Y_0} \in \Sigma$. The fold-saddle singularity $p_0=S_{Y_0}$
is illustrated in Figures \ref{fig fold-sela} and \ref{fig
fold-visivel-sela} $-$ the dotted lines in these and later figures
represent the points where $X.f=0$ and $Y.f=0$.

Let $p=(0,0)$ be a fold-saddle singularity of $Z=(X,Y)$. 
We denote the set of all NSVF $Z=(X,Y)$ such
that $X \in \Gamma^{F}_{\Sigma_{+}}$  
and $Y \in \Gamma^{S}_{\Sigma_{-}}$ by $\Gamma^{F-S}$.
We endow $\chi^r \times \chi^r$ (consequently $\Omega^r$ and
$\Gamma^{F-S}$) with the product topology. Let
$Z_0=(X_0,Y_0) \in \Gamma^{F-S}$. 
Observe that $0$ is the unique singularity of $X_0$ around a
neighborhood $W_0$ of the origin in $\R^2$. So, there exists a
neighborhood $U_0$ of $Z_0$ in $\Omega^r$ such that for any $Z=(X,Y)
\in U_0$ we may find a $\Sigma-$fold point $p_Z=(k_Z,0) \in W_0$
such that it is the unique singularity of $X$ in $W_0.$ Moreover the
correspondence $Z \mapsto p_Z$ is $C^r.$

In the same way, for any $Z=(X,Y)\in U_0$ we find a
$C^r-$correspondence $B:U_0\rightarrow \R^2$ where
$B(Z)=s_Z=(a_Z,b_Z)$ is the (unique) equilibrium (saddle) of $Y$ in
$U_0$. We are assuming that the eigenspaces of $DY_{s_Z} (q_Z)$ are
transverse to $\Sigma$ at $s_Z$. We have to distinguish the cases:
$(i)$ $b_Z<0$, $(ii)$ $b_Z=0$ and $(iii)$ $b_Z>0.$ Observe that when
$b_Z<0$ (resp. $b_Z>0$) there is associated to $Z$ an invisible
(resp. visible) $\Sigma-$fold point of $Y$ given by $q_Y=(c_Z,0) \in
W_0$. Moreover $lim_{b_Z\rightarrow 0}c_Z = a_Z.$

Define $F(Z)=(k_Z-a_Z,b_Z)$. Knowing that $s_Z$ is a hyperbolic
equilibrium point of $Z$ and $0$ is a regular value of $F$, it is
not hard to prove (see \cite{Marco-bollettino}) that:
\begin{itemize}
\item The derivative $DF:U_0 \rightarrow \R^2$ is surjective
and

\item $F^{-1}(0)=\Omega_2$ is a codimension two submanifold of
$\Omega^r$.\end{itemize} Therefore  this fold-saddle singularity
occurs generically in
two-parameter fa\-mi\-lies of vector fields in $\Omega^r$. 

\subsection{Normal Form}\label{subsecao forma normal}

We start this section with the following model:
\begin{equation}\label{eq fold-sela zero}
Z^{\tau} = \left\{
      \begin{array}{ll}
        X^{\tau} = \left(
              \begin{array}{c}
                   \rho_1  \\
                \alpha_1(\tau)x
\end{array}
      \right)
 & \hbox{if $y \geq 0$,} \\
         Y = \left(
              \begin{array}{c}
               k_1 y \\
               k_1 x
\end{array}
      \right)& \hbox{if $y \leq 0,$}
      \end{array}
    \right.
\end{equation}
where $\tau= \{ inv , vis \}$, $\alpha_1(inv)=-1$,
$\alpha_1(vis)=1$, $\rho_1=  \, \pm 1$ and $k_1 = \, \pm 1$.\\

The next lemma provides explicitly the equivalence between $Z \in
\Omega_2$ and the model (\ref{eq fold-sela zero}). We present an
outline of proof of  the previous lemma in Section \ref{secao prova
lema 1}.

\begin{lemma}\label{lema equivalencia} If $Z \in \Omega_2$ then $Z$ is $\Sigma-$equivalent to
$Z^{\tau}$ given by  (\ref{eq fold-sela zero}).
\end{lemma}

%

Note that there exists a NSVF $\widetilde{Z} \in \Omega^r$ nearby
$Z^{inv}$, given by  (\ref{eq fold-sela zero}),  such that
$\widetilde{Z}$ presents a $\Sigma-$center (see Figure \ref{fig
sigma centro}). In fact, this suggests that the unfolding of
(\ref{eq fold-sela zero}) has infinite codimension. At this point
it seems natural to propose the following conjecture.\\

\noindent{\textbf{Conjecture: }}\textit{ For any neighborhood
$\mathcal{W} \subset \Omega^r$ of $Z^{inv}$ (given by  (\ref{eq
fold-sela zero})), and for any integer $k>0$ there exists
$\widetilde{Z} \in W$ such
that the codimension of $\widetilde{Z}$ is $k$.}\\

So, based on this conjecture, we have to sharpen our normal normal.
In fact, in order to get low codimension bifurcation we have to
impose some
generic assumptions.\\

Without loss of generality, throughout the rest of this paper we
consider $\rho_1 =1$ and $k_1=-1$ at the model (\ref{eq fold-sela
zero}). The other choices of the parameters are treated analogously.  
When $\tau=inv$ we add the extra generic assumption $X_{0}^{3}.f(p)
\neq
0$ on the $\Sigma-$fold point $p$ of $Z_0=(X_0,Y_0) \in \Gamma^{F-S}$. 
By means of Theorem $2$ in \cite{vishik}, 
 we may conclude that around  the invisible
$\Sigma-$fold point the vector field $X_0$ can be expressed as
$X_0=(1, - x + a_1 x^2)$ and $f(x,y)=y$, where $a_1\neq0$. We say
that  $X_0$ is \textit{contractive} (respectively,
\textit{expansive}) at $p$ if
$a_1<0$ (respectively $a_1>0$).\\

According to the previous discussion, we will consider $Z^{inv}_0,
Z^{vis}_0 \in \Omega^r$ written in the following forms:
\begin{equation}\label{eq fold-sela inicio}
Z^{inv}_0 = \left\{
      \begin{array}{ll}
        X^{inv}_0 = \left(
              \begin{array}{c}
                1 \\
               -x + x^2
\end{array}
      \right)
 & \hbox{if $y \geq 0$,} \\
        Y_0 =  \left(
              \begin{array}{c}
                -y \\
               -x
\end{array}
      \right)& \hbox{if $y \leq 0$, and}
      \end{array}
    \right.
\end{equation}
\begin{equation}\label{eq fold-visivel-sela inicio}
Z^{vis}_0 = \left\{
      \begin{array}{ll}
        X^{vis}_0 = \left(
              \begin{array}{c}
                1 \\
               x
\end{array}
      \right)
 & \hbox{if $y \geq 0$,} \\
        Y_0 =  \left(
              \begin{array}{c}
                -y \\
               -x
\end{array}
      \right)& \hbox{if $y \leq 0$.}
      \end{array}
    \right.
\end{equation}
Note that $X^{inv}_0$ presents an invisible expansive $\Sigma-$fold
point in its phase portrait and $X^{vis}_0$ presents a visible one.

\subsection{Unfoldings}\label{subsecao desdobramentos}

The main question of this paper is to exhibit the bifurcation
diagram of
$Z^{\tau}_0$ with either $\tau=inv$ or $\tau=vis$. 
%
For this reason we consider unfoldings of the normal forms (\ref{eq
fold-sela inicio}) and (\ref{eq fold-visivel-sela inicio}). We
obtain that:

I- There is an imbedding $F_{0}^{\tau}:\R^2,0 \rightarrow
\chi^r,Z^{\tau}_0$ such that
$F_{0}^{\tau}(\lambda,\beta)=Z^{\tau}_{\lambda,\beta}$ is expressed
by:
\begin{equation}\label{eq fold-sela sem variar}
Z^{\tau}_{\lambda, \beta} = \left\{
      \begin{array}{ll}
        X^{\tau}_{\lambda} = \left(
              \begin{array}{c}
                   1 \\
                \alpha_1(\tau)(x - \lambda) + \alpha_2(\tau)(x-\lambda)^2
\end{array}
      \right)
 & \hbox{if $y \geq 0$,} \\
         Y_{\beta} = \left(
              \begin{array}{c}
               -(y + \beta) \\
               - x
\end{array}
      \right)& \hbox{if $y \leq 0,$}
      \end{array}
    \right.
\end{equation}
where  $\tau= \{ inv , vis \}$, $(\lambda, \beta) \in (-1,1) \times
(-\sqrt{3}/2,\sqrt{3}/2)$, $\alpha_1(inv)=-1$, $\alpha_1(vis)=1$,
$\alpha_2(inv)=1$ and $\alpha_2(vis)=0$. Moreover, the two-parameter
family given by
 (\ref{eq fold-sela sem variar}) is transversal to $\Omega_2$. 
 We stress that in
\cite{Marcel}  the unfolding of the case $\tau=inv$ is done.
Nevertheless we observe that there are some typical topological
types nearby $Z^{\tau}_0$ that do not appear in the bifurcation
diagram of $Z^{\tau}_{\lambda,\beta}$. For exam\-ple, when
$\tau=inv$ the configurations in Fi\-gu\-re \ref{fig nao aparecem}
cases (a) and (b) are excluded and when $\tau=vis$ the configuration
in Figure \ref{fig nao aparecem}$-(c)$ also is excluded.

\begin{figure}[h]
\begin{center}
\psfrag{a}{(a)} \psfrag{b}{(b)} \psfrag{c}{(c)} \psfrag{D}{}
\psfrag{E}{\hspace{1cm}$Z^\Sigma(q)$}
\psfrag{F}{\hspace{-.35cm}$\Sigma_s$} \psfrag{G}{} \epsfxsize=15cm
\epsfbox{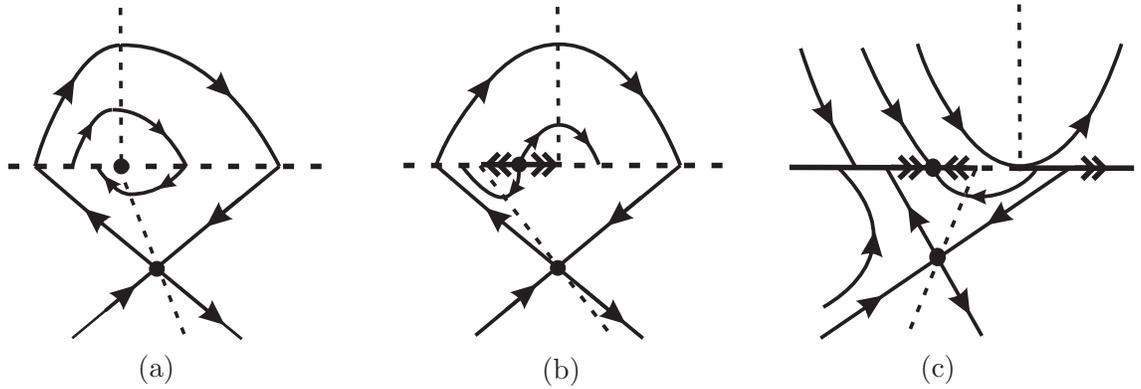} \caption{\small{Cases do not covered
in \cite{Marcel}. In (a) and (b) (respectively, in (c)) the fold
singularity is invisible (respectively, visible).}} \label{fig nao
aparecem}
\end{center}
\end{figure}

%

II- In order to consider a more general situation which includes the
above mentioned cases we add an auxiliary parameter $\mu$. As a
result the model in Equation \ref{eq fold-sela com parametros} is
obtained.\\


We stress that the configuration illustrated in Figure \ref{fig nao
aparecem}$-(a)$ plays a very important role in our analysis. In this
\textit{resonant} configuration we note, simultaneously, a fold-fold
singularity (which is a local phenomenon) and a $\Sigma-$graph
(loop) passing through the saddle equilibrium (which is a global
phenomenon). Only the bifurcation of these two unstable
configurations already represents a relevant development (the
fold-fold singularity was studied recently in \cite{Marcel} and the
\textit{non-smooth loop bifurcation}, as far as we know, was not
studied until the present work). In fact, this configuration is
reached in  (\ref{eq fold-sela com parametros}), taking
$\mu=\mu_0(\beta)$.\\ 

In what follows, in order to simplify the calculations, we take $\mu
=\alpha + 1$ in  (\ref{eq fold-sela com parametros}) and obtain the
following expression:
\begin{equation}\label{eq fold-sela com parametros novos}
Z^{\tau}_{\lambda , \alpha, \beta} = \left\{
      \begin{array}{ll}
        X_{\lambda} = \left(
              \begin{array}{c}
                  1 \\
                \alpha_1(\tau)(x - \lambda) + \alpha_2(\tau)
                (x-\lambda)^2
\end{array}
      \right)
 & \hbox{if $y \geq 0$,} \\
         Y_{\alpha , \beta} = \left(
              \begin{array}{c}
                \frac{(1 + \alpha)}{2} x + \frac{(-1 + \alpha)}{2}(y + \beta) \\
                \frac{(-1 + \alpha)}{2} x + \frac{(1 + \alpha)}{2}(y + \beta)
\end{array}
      \right)& \hbox{if $y \leq 0$,}
      \end{array}
    \right.
\end{equation}
where  $\tau= \{ inv, vis \}$, $\alpha_1(inv)=-1$,
$\alpha_1(vis)=1$, $\alpha_2(inv)=1$, $\alpha_2(vis)=0$ and
$(\lambda,\beta,\alpha) \in (-1,1) \times (-\sqrt{3}/2,\sqrt{3}/2)
\times (-1 - \varepsilon_0, 1)$ with $\varepsilon_0 > 0$
sufficiently small. Since $\mu_0(\beta)$ is given by  (\ref{eq mu
zero}), we obtain that
\begin{equation}\label{eq alpha zero}\alpha_0(\beta) = 1 - (12 \beta/(-3 +
6 \beta + \sqrt{9- 12 \beta^2})).\end{equation}

When it does not produce confusion, in order to simplify the
notation we
 use  $Z=(X,Y)$ or $Z_{\lambda , \alpha,
\beta} = (X , Y)$ instead $Z^{\tau}_{\lambda , \alpha, \beta} =
(X_{\lambda}, Y_{\alpha , \beta})$.

\subsubsection{Geometrical Analysis of the Normal Form (\ref{eq fold-sela com parametros
novos})}\label{subsecao analise geometrica forma normal}

Given $Z=(X,Y)$, we describe some properties of both $X =
X_{\lambda}$ and $Y = Y_{\alpha, \beta}$.

The real number $\lambda$ measures how the $\Sigma-$fold point
$d=(d_1,d_2)=(\lambda,0)$ of $X$ is translated away from the origin.
More specifically, if $\lambda < 0$ then $d$ is translated to the
left hand side and if $\lambda
> 0$ then $d$ is
translated to the right hand side.

\begin{figure}[!h]
\begin{minipage}[b]{0.485\linewidth}
\begin{center}\psfrag{A}{$\Sigma$} \psfrag{B}{$e$} \psfrag{G}{$Y.f=0$}
\psfrag{D}{$\Sigma$} \psfrag{E}{$X.f=0$} \psfrag{F}{$\gamma_1$}
\epsfxsize=4cm \epsfbox{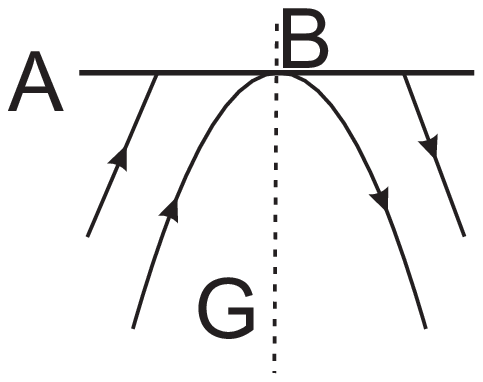} \caption{\small{Case
$Y^-$: perturbation of a boundary saddle generating a visible
fold.}} \label{fig sela para cima}\end{center}
\end{minipage} \hfill
\begin{minipage}[b]{0.485\linewidth}
\begin{center}\psfrag{A}{$h$} \psfrag{B}{$i$} \psfrag{C}{$j$}
\psfrag{D}{$\Sigma_e$} \psfrag{E}{$\Sigma_c$} \psfrag{F}{$S$}
\psfrag{G}{$Y.f=0$} \epsfxsize=3cm \epsfbox{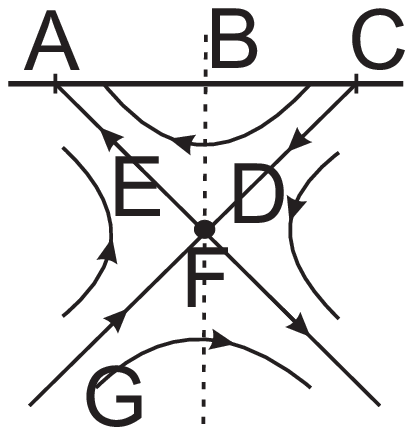}
\caption{\small{Case $Y^+$: perturbation of a boundary saddle
generating an invisible fold.}} \label{fig sela para baixo}
\end{center}\end{minipage} \hfill
\end{figure}

Some calculations show that the curve $Y.f=0$ is given by
$y=\frac{(1 - \alpha)}{(1 + \alpha)}x - \beta$. So the points of
this curve are equidistant from the separatrices when $\alpha = -1$.
It becomes closer to the stable separatrix of the saddle equilibrium
$S = S_{\alpha, \beta}=(s_1,s_2)$ when $\alpha \in (-1,0)$. It
becomes closer to the unstable separatrix of $S$ when $\alpha \in
(-1 +\varepsilon_0,-1)$. Moreover, the smooth vector field $Y$ has
distinct types of contact with $\Sigma$ according to the particular
deformation considered. In this way, we have to consider the
following behaviors:

\begin{itemize}
\item $\mathbf{Y^- :}$ In this case $\beta < 0$. So  $S$ is translated to the
$y-$direction with $y>0$ (and $S$ is not visible for $Z$). It has a
visible $\Sigma-$fold point $e= e_{\alpha, \beta} = (e_1,e_2)=\Big(
\frac{(1+ \alpha)}{(1 - \alpha)}\beta,0 \Big)=(e_1,0)$ (see Figure
\ref{fig sela para cima}).

\item $\mathbf{Y^0 :}$ In this case $\beta = 0$. So $S$ is not translated (see Figure \ref{fig fold-sela}).

\item $\mathbf{Y^+ :}$ In this case $\beta > 0$. So $S$ is translated to the
$y-$direction with $y<0$. It has an invisible $\Sigma-$fold point $i
=  i_{\alpha, \beta} = (i_1,i_2) = \Big( \frac{(1+ \alpha)}{(1 -
\alpha)}\beta,0 \Big)$. Moreover, we distinguish two  points:  $h =
h_{\beta}= (h_1,h_2) \linebreak = (-\beta,0)$ which is the
intersection between the unstable separatrix with $\Sigma$ and  $j =
j_{\beta}= (j_1,j_2)=(\beta,0)$ which is the intersection between
the stable separatrix with $\Sigma$ (see Figure \ref{fig sela para
baixo}).
\end{itemize}

 In Figure \ref{fig sela para baixo} we
distinguish the arcs of trajectory $\sigma_{1}$ joining  $S$ to $h$
and  $\sigma_{2}$  joining $j$ to $S$.

\section{Proof of Lemma 1}\label{secao prova lema 1}


Now, we show how we can construct the homeomorphism in Lemma
\ref{lema equivalencia}.\\

\begin{proof}[Outline of Proof of Lemma 1]
Here we construct a $\Sigma-$preserving homeomorphism $h$ that sends
orbits of $Z=(X,Y) \in \Omega_2$, defined in a sufficiently small
neighborhood $\mathcal{U}_0$ of the fold-cusp singularity of $Z$, to
orbits of $\widetilde{Z}=(\widetilde{X},\widetilde{Y})$, defined in
a sufficiently small neighborhood $\widetilde{\mathcal{U}}_0$ of the
fold-cusp singularity of $\widetilde{Z}$, where
$\widetilde{Z}=Z^{inv}$ is given by  (\ref{eq fold-sela zero}) with
$\rho_1 = 1$ and $k_1= -1$ (the other cases are treated
analogously). Consider $A_0$ an arbitrary point of the stable
separatrix of the saddle point $S$ of $Y$ (see Figure \ref{fig
homeomorfismo}). Let $T_1$ be a transversal section of $Y$ at $A_0$.
The section $T_1$ also is transversal to $\widetilde{Y}$ and it
crosses the stable separatrix of the saddle point $\widetilde{S}$ of
$\widetilde{Y}$ at $B_0$. Let $A_1 \in T_1$ be a point on the left
of $A_0$. The trajectory of $Y$ passing through $A_1$ crosses
$\Sigma$ at $A_2$. In the same way, the trajectory of
$\widetilde{Y}$ passing through $B_1$ crosses $\Sigma$ at $B_2$. The
trajectory of $X$ that passes through $A_2$ crosses $\Sigma$ in a
point $A_3$.
 The trajectory of $\widetilde{X}$ that passes through $B_2$ crosses
 $\Sigma$ in $B_3$. Consider $A_4$ an arbitrary
 point of the unstable separatrix of $S$. Let $T_2$ be a transversal section of $Y$ passing through
 $A_4$. The section $T_2$ also is transversal to $\widetilde{Y}$ and
 it crosses the unstable separatrix of $\widetilde{S}$ at the point
 $B_4$. The
 trajectory of $Y$ passing through $A_3$ crosses $T_2$ in a point
 $A_5$. In the same way, the trajectory of $\widetilde{Y}$ passing
 through $B_3$ crosses $T_2$ at $B_5$. Let $A_6 \in T_1$ be a point at the right of $A_0$. The
 trajectory of $Y$ passing through $A_6$ crosses $T_2$ at $A_7$. The trajectory of $\widetilde{Y}$ passing
 through $A_6$ crosses $T_2$ at $B_7$. The homeomorphism $h$
 sends $T_1$ to $T_1$, the arc of trajectory $\gamma_1=\widehat{A_1 \, A_5}$
 to the arc of trajectory $\widetilde{\gamma}_1=\widehat{A_1 \, B_5}$ and the arc of trajectory $\gamma_2=\widehat{A_6 \, A_7}$
 to the arc of trajectory $\widetilde{\gamma}_2=\widehat{A_6 \, B_7}$. Now we can extend
 continuously $h$ to the interior of the region bounded by
 $T_1 \cup \gamma_1 \cup T_2 \cup \gamma_2$. In this way, there exists a $\Sigma-$preserving homeomorphism $h$
 that sends orbits of $Z$ to orbits of $\widetilde{Z}$.

\begin{figure}[!h]
\begin{center}
\psfrag{A}{$A_1$} \psfrag{B}{$A_6$} \psfrag{C}{$A_0$}
\psfrag{D}{$T_1$} \psfrag{E}{$A_2$} \psfrag{F}{$A_3$}
\psfrag{G}{$A_5$} \psfrag{H}{$A_7$} \psfrag{I}{$A_4$}
\psfrag{J}{$T_2$} \psfrag{K}{$B_0$} \psfrag{L}{$B_2$}
\psfrag{M}{$B_3$} \psfrag{N}{$B_5$} \psfrag{O}{$B_4$}
\psfrag{P}{$B_7$} \psfrag{Q}{$h$}  \epsfxsize=12cm
\epsfbox{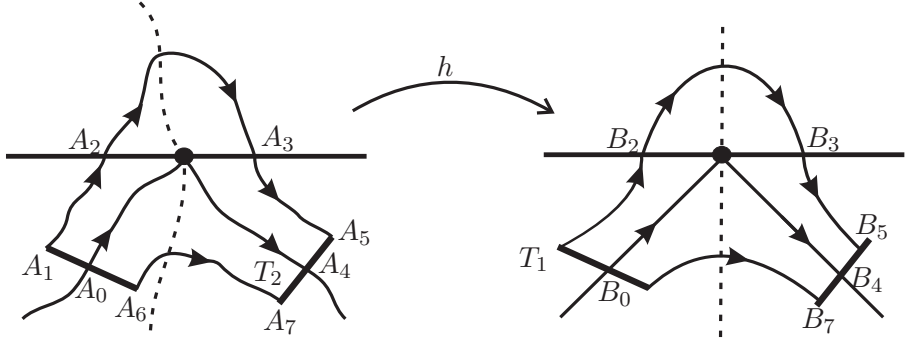} \caption{\small{Construction of a
$C^0-$equivalence: Case $\tau=inv$.}} \label{fig homeomorfismo}
\end{center}
\end{figure}

\noindent When  $\widetilde{Z}=Z^{vis}$ is given by  (\ref{eq
fold-sela zero}), the first coordinate of $\widetilde{X}$ is equal
to $1$ and $k_1= -1$. We build the same above construction until the
appearance of $A_2$ and $B_2$. Now, consider $C_4$ an arbitrary
point of the stable separatrix of the $\Sigma-$fold point $F$ of $X$
(see Figure \ref{fig homeomorfismo visivel}). Let $T_3$ be a
transversal section to $X$ at $C_4$. The section $T_3$ is also
transversal to $\widetilde{X}$ and it crosses the stable separatrix
of the $\Sigma-$fold point $\widetilde{F}$ of $\widetilde{X}$ at the
point $D_4$. The trajectory of $X$ passing through $A_2$ crosses
$T_3$ at a point $C_3$. In the same way, the trajectory of
$\widetilde{X}$ passing through $B_2$ crosses $T_3$ at a point
$D_3$. Consider $C_5$ an arbitrary point of the unstable separatrix
of $F$. Let $T_4$ be a transversal section to $X$ at $C_5$. The
section $T_4$ also is transversal to $\widetilde{X}$ and it crosses
the unstable separatrix of $\widetilde{F}$  at  $D_5$. Let $C_6 \in
T_3$ be a point at the right of $C_4$. The trajectory of $X$ passing
through $C_6$ crosses $T_4$ at $C_7$. In the same way, let $D_6 \in
T_3$ be a point at the right of $D_4$. The trajectory of
$\widetilde{X}$ passing through $D_6$ crosses $T_4$ at $D_7$.  Let
$C_8 \in T_4$ be a point at the right of $C_5$. The trajectory of
$X$ passing through $C_8$ crosses $\Sigma$ at $A_3$. In the same
way, let $D_8 \in T_4$ be a point at the right of $D_5$. The
trajectory of $\widetilde{X}$ passing through $D_8$ crosses $\Sigma$
at $B_3$. Now, it is enough to repeat the construction made in the
previous case. The homeomorphism $h$ sends $T_1$ to $T_1$, the arc
of trajectory $\gamma_1=\widehat{A_1 \, C_3}$ to the arc of
trajectory $\widetilde{\gamma}_1=\widehat{A_1 \, D_3}$, $T_3$ to
$T_3$, the arc of trajectory $\gamma_2=\widehat{C_6 \, C_7}$ to the
arc of trajectory $\widetilde{\gamma}_2=\widehat{D_6 \, D_7}$, $T_4$
to $T_4$, the arc of trajectory $\gamma_3=\widehat{C_8 \, A_5}$ to
the arc of trajectory $\widetilde{\gamma}_3=\widehat{D_8 \, B_5}$,
$T_2$ to $T_2$ and the arc of trajectory $\gamma_4=\widehat{A_6 \,
A_7}$ to the arc of trajectory $\widetilde{\gamma}_4=\widehat{A_6 \,
B_7}$. Now we can extend continuously $h$ to the interior of the
region bounded by $T_1 \cup \gamma_1 \cup T_3 \cup \gamma_2 \cup T_4
\cup \gamma_3 \cup T_2 \cup \gamma_4$. In this way, there exists a
$\Sigma-$preserving homeomorphism $h$ that sends orbits of $Z$ to
orbits of $\widetilde{Z}$.\end{proof}

\begin{figure}[!h]
\begin{center}
\psfrag{A}{$A_1$} \psfrag{B}{$A_6$} \psfrag{C}{$A_0$}
\psfrag{D}{$T_1$} \psfrag{E}{$A_2$} \psfrag{F}{$A_3$}
\psfrag{G}{$A_5$} \psfrag{H}{$A_7$} \psfrag{I}{$A_4$}
\psfrag{J}{$T_2$} \psfrag{K}{$B_0$} \psfrag{L}{$B_2$}
\psfrag{M}{$B_3$} \psfrag{N}{$B_5$} \psfrag{O}{$B_4$}
\psfrag{P}{$B_7$} \psfrag{Q}{$h$} \psfrag{u}{$T_3$}
\psfrag{r}{$C_3$} \psfrag{s}{$C_4$} \psfrag{t}{$C_6$}
\psfrag{2}{$T_4$} \psfrag{v}{$C_7$} \psfrag{0}{$C_5$}
\psfrag{1}{$C_8$} \psfrag{3}{$D_3$} \psfrag{4}{$D_4$}
\psfrag{5}{$D_6$} \psfrag{6}{$D_7$} \psfrag{7}{$D_5$}
\psfrag{8}{$D_8$}  \epsfxsize=11cm
\epsfbox{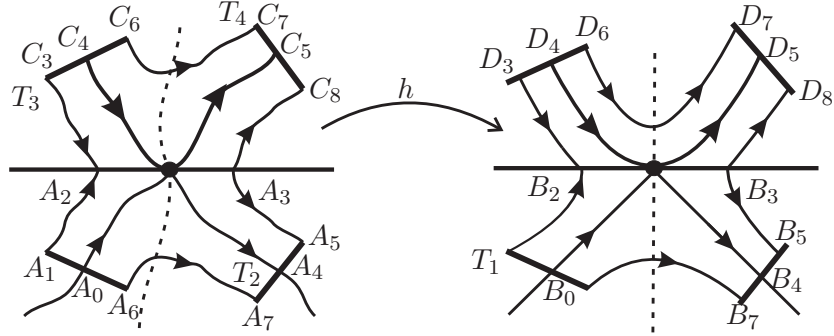} \caption{\small{Construction
of a $C^0-$equivalence: Case $\tau=vis$.}} \label{fig homeomorfismo
visivel}
\end{center}
\end{figure}


\section{Proof of Theorem 1}\label{secao prova teorema 1}

\begin{proof}[Proof of Theorem 1] In Cases $1_1$, $2_1$ and $3_1$ we assume
that $Y$ presents the behavior $Y^-$. In Cases $4_1$, $5_1$ and
$6_1$ we assume that $Y$ presents the behavior $Y^0$. In these cases
canard cycles do not arise.

$\diamond$ \textit{Cases $(1_1)$ $d_1<e_1$, $(2_1)$ $d_1=e_1$ and
$(3_1)$ $d_1 > e_1$:} The points of $\Sigma$ outside the interval
$(d_1,e_1)$ (or  $(e_1 ,d_1)$) belong to $\Sigma_c$. The points
inside this interval, when it is not degenerate, belong to
$\Sigma_s$ in Case $1_1$ and to $\Sigma_{2}$ in Case $3_1$. In both
cases $H(z)>0$ for all $z \in \Sigma_e \cup \Sigma_s$. See Figure
\ref{fig 1 teo 1}.

\begin{figure}[!h]
\begin{center}\psfrag{A}{\hspace{-.5cm}$\lambda < 0$} \psfrag{B}{$\lambda=0$}
\psfrag{C}{$\lambda > 0$}\psfrag{0}{$1_1$} \psfrag{1}{$2_1$}
\psfrag{2}{$3_1$} \epsfxsize=9.5cm \epsfbox{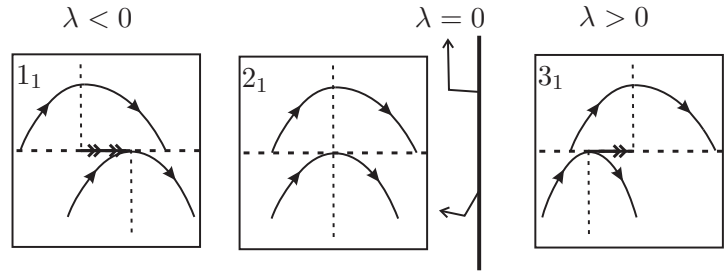}
\caption{\small{Cases $1_1$, $2_1$ and $3_1$.}} \label{fig 1 teo 1}
\end{center}
\end{figure}

$\diamond$ \textit{Cases $(4_1)$ $d_1 < s_1$, $(5_1)$ $d_1=s_1$ and
$(6_1)$ $d_1 > s_1$:} The points of $\Sigma$ outside the interval
$(d_1,s_1)$ (or $(s_1 ,d_1)$) belong to $\Sigma_c$. The points
inside this interval, when it is not degenerate, belong to
$\Sigma_s$ in Case $4_1$ and to $\Sigma_{2}$ in Case $6_1$. In both
cases $H(z)>0$ for all $z \in \Sigma_e \cup \Sigma_s$. See Figure
\ref{fig 2 teo 1}.

\begin{figure}[!h]
\begin{center}\psfrag{A}{$\lambda < 0$} \psfrag{B}{$\lambda=0$} \psfrag{C}{$\lambda > 0$}\psfrag{0}{$4_1$} \psfrag{1}{$5_1$} \psfrag{2}{$6_1$}
\epsfxsize=9.5cm \epsfbox{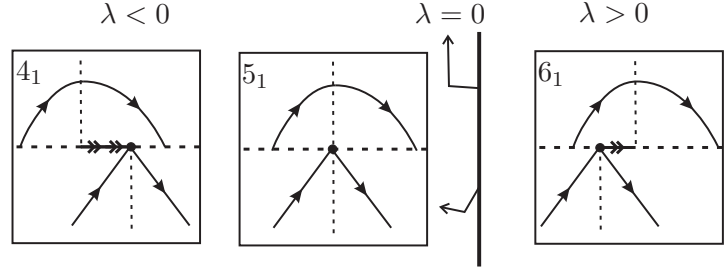} \caption{\small{Cases
$4_1$, $5_1$ and $6_1$.}} \label{fig 2 teo 1}
\end{center}
\end{figure}

In Cases $7_1 - 19_1$ we assume that $Y$ presents the behavior
$Y^+$.

Observe that as the parameter $\lambda$ increases $-$ assuming the
values $L_0$, $L_1$ and $L_2$ below described $-$ it appear
orbit-arcs of $X$ connecting the points $h$ and $i$, $h$ and $j$ and
$i$ and $j$ respectively.


Remembering that $\alpha =  1 - (12 \beta/(-3 + 6 \beta + \sqrt{9-
12 \beta^2}))$ the values of $L_0$, $L_1$ and $L_2$ are:

$\begin{array}{ccl}
  L_0  & = &
[-9 -6 \beta + \sqrt{9-12 \beta^2} + \sqrt{2} \sqrt{15 + \sqrt{9-12
\beta^2} - 2 \beta (-3 + 2 \beta + \sqrt{9 - 12 \beta^2})}]/12 \\
   L_1 & = & -1/2 +
\sqrt{9 - 12 \beta^2}/6 \\
   L_2 & = & [-9 +6 \beta + \sqrt{9-12 \beta^2} +  \sqrt{2} \sqrt{15 +
\sqrt{9-12 \beta^2} + 2 \beta (-3 + 2 \beta + \sqrt{9 - 12
\beta^2})}]/12. \\
 \end{array}$

$\diamond$ \textit{Cases $(7_1)$ $\lambda< - \beta$, $(8_1)$
$\lambda = - \beta$, $(9_1)$ $- \beta < \lambda < L_0$, $(10_1)$
$\lambda = L_0$ and $(11_1)$ $L_0 < \lambda < L_1 $:} The points of
$\Sigma$ outside the interval $(d_1,i_1)$ belong to $\Sigma_c$. The
points inside this interval belong to $\Sigma_s$. The direction
function $H$ assumes positive values in a neighborhood of $d_1$,
negative values in a neighborhood of $i_1$ and there exists only one
value $\widetilde{P}=\widetilde{P}_{\lambda,\alpha,\beta}$ such that
$H(\widetilde{P})=0$. 
So, by  (\ref{eq H}), the $\Sigma-$attractor $P=(\widetilde{P},0)$,
nearby $(0,0)$, is the unique pseudo equilibrium of $Z$. In these
cases canard cycles do not arise. See Figure \ref{fig 3a teo 1}.

\begin{figure}[!h]
\begin{center}\psfrag{A}{\hspace{-.8cm}$\lambda < - \beta$} \psfrag{B}{\hspace{-.4cm}$\lambda= - \beta$}
\psfrag{C}{\hspace{-.35cm}$- \beta < \lambda < L_0$} \psfrag{D}{$
\lambda = L_0$} \psfrag{E}{\hspace{.35cm}$L_0 < \lambda <
L_1$}\psfrag{0}{$7_1$} \psfrag{1}{$8_1$} \psfrag{2}{$9_1$}
\psfrag{3}{$10_1$}\psfrag{4}{$11_1$} \epsfxsize=10cm
\epsfbox{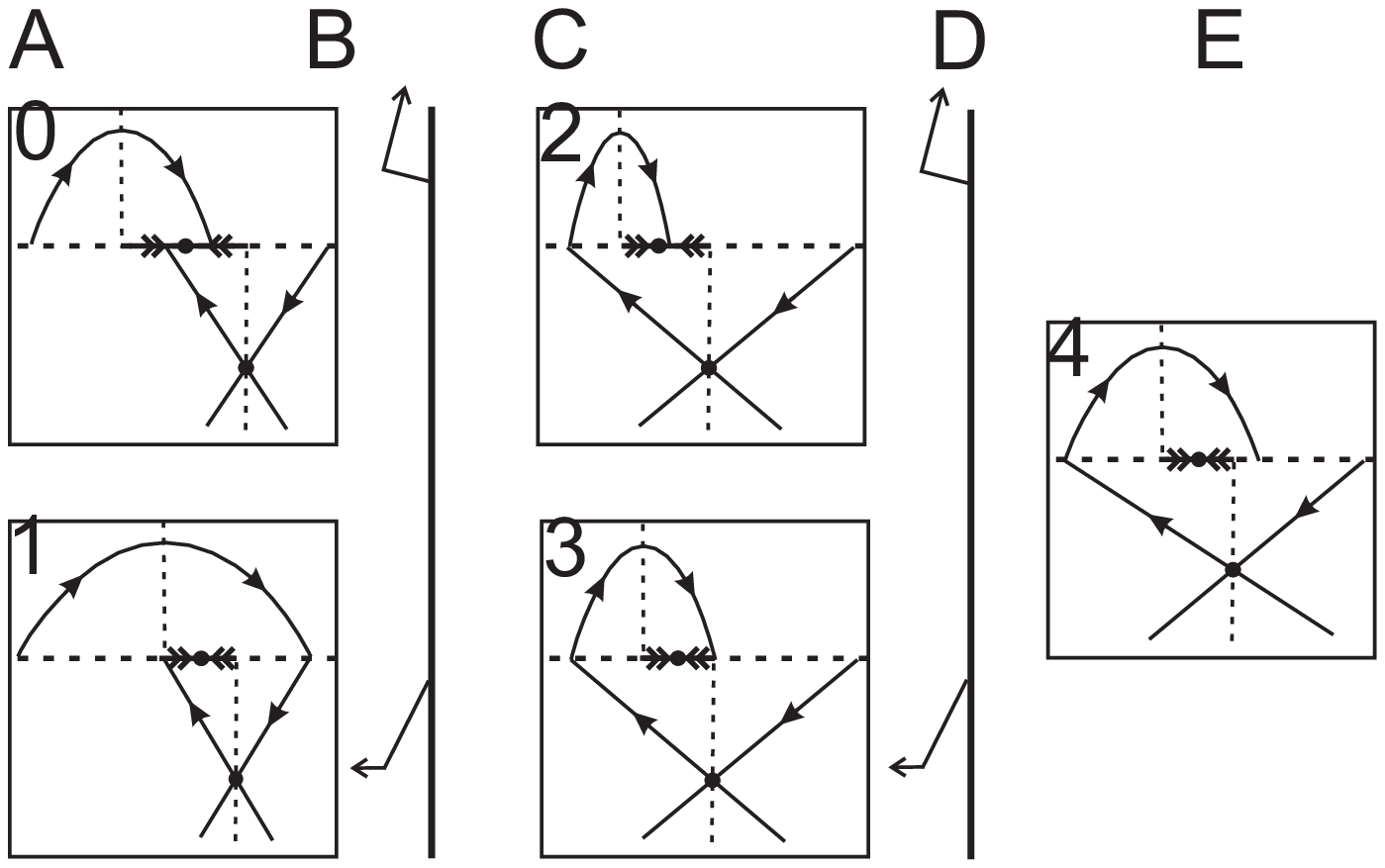} \caption{\small{Cases $7_1 - 11_1$.}}
\label{fig 3a teo 1}
\end{center}
\end{figure}

$\diamond$ \textit{Case $(12_1)$ $\lambda = L_1$:} Since
 $\lambda=L_1$ there is an orbit-arc
$\gamma_{1}^X$ of $X$ connecting the points $h$ and $j$. It
generates a $\Sigma-$graph $\Gamma = \gamma_{1}^{X} \cup \Sigma_e
\cup S \cup \Sigma_c$ of kind I. Moreover, since $\alpha=\alpha_0$,
where $\alpha_0$ is given by (\ref{eq alpha zero}), there exists a
non generic tangential singularity at the point $d=i$. So, the
points of $\Sigma/\{d \}$ belong to $\Sigma_c$. As the $\Sigma-$fold
point of $X$ is expansive, a direct calculus shows that the
\textit{First Return Map} $\eta: (h,d) \rightarrow (h,d)$ has
derivative bigger than $1$. As consequence, $\Gamma$ is a repeller
for the trajectories inside it, $d=i$ behaves itself like an
attractor (weak focus) and canard cycles do not arise.  See Figure
\ref{fig nova teo 1}.


\begin{figure}[!h]
\begin{center}\psfrag{A}{$\lambda =L_1$} \psfrag{B}{$L_1<\lambda < L_3$} \psfrag{C}{$\lambda = L_3$} \psfrag{D}{$ \lambda = -\frac{\beta}{2}$}
\psfrag{E}{$-\frac{\beta}{2} < \lambda < 0$}\psfrag{0}{$12_1$}
\psfrag{1}{$13_1$} \psfrag{2}{$14_1$}
\psfrag{3}{$10_1$}\psfrag{4}{$11_1$} \epsfxsize=12.5cm
\epsfbox{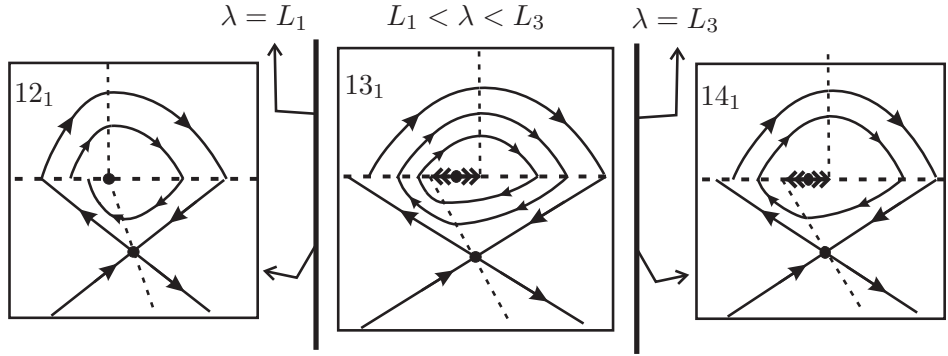} \caption{\small{Cases $12_1$, $13_1$
and $14_1$.}} \label{fig nova teo 1}
\end{center}
\end{figure}

$\diamond$ \textit{Case $(13_1)$ $L_1< \lambda < L_3$:} The meaning
of  $L_3$ will be given below in this case. The points of $\Sigma$
outside the interval $(i_1,d_1)$ belong to $\Sigma_c$ and the points
inside this interval belong to $\Sigma_e$. The direction function
$H$ assumes positive values in a neighborhood of $d_1$, negative
values  in a  neighborhood of $i_1$ and there exists a unique value
$\widetilde{P}=\widetilde{P}_{\lambda,\alpha,\beta}$ such that
$H(\widetilde{P})=0$. So $P=(\widetilde{P},0)$ is a
$\Sigma-$repeller. When $\lambda$ is a bit bigger than $L_1$, the
First Return Map $\eta$ has two fixed points, i.e., $Z$ has two
canard cycles. One of them, called $\Gamma_1$, is born from the
bifurcation of the $\Sigma-$graph $\Gamma$ of the previous case and
the other one, called $\Gamma_2$, is born from the bifurcation of
the non generic tangential singularity presented in the previous
case. Both of them are canard cycles of kind I. Moreover, we obtain
that $\Gamma_1$ is a hyperbolic repeller canard cycle and $\Gamma_2$
is a hyperbolic attractor canard cycle. Note that, as $\lambda$
increases, $\Gamma_1$ becomes smaller and $\Gamma_2$ becomes bigger.
When $\lambda$ assumes the limit value $L_3$, one of
them collides to the other. 
See Figure \ref{fig nova teo 1}.

$\diamond$ \textit{Case $(14_1)$ $\lambda = L_3$:} The distribution
of the connected components of $\Sigma$ and the behavior of $H$ are
the same as Case $13_1$. Since $\lambda=L_3$, as described in the
previous case, there exists a non hyperbolic canard cycle $\Gamma$
of kind I which is an attractor for the trajectories inside it and
is a repeller for the trajectories outside it.  See Figure \ref{fig
nova teo 1}.

$\diamond$ \textit{Cases $(15_1)$ $L_3<\lambda< L_2$, $(16_1)$
$\lambda = L_2$, $(17_1)$ $L_2 < \lambda < \beta$, $(18_1)$ $\lambda
= \beta$ and $(19_1)$ $\lambda
> \beta$:} The points of $\Sigma$ outside the interval $(i_1,d_1)$
belong to $\Sigma_c$ and the points inside this interval belong to
$\Sigma_e$. The direction function $H$ assumes positive values in a
neighborhood of $d_1$, negative values in a neighborhood of $i_1$
and there exists a unique value $\widetilde{P}$ such that
$H(\widetilde{P})=0$. So, by  (\ref{eq H}), the $\Sigma-$repeller
$P=(\widetilde{P},0)$, nearby $(0,0)$, is the unique pseudo
equilibrium of $Z$. In these cases canard cycles do not arise. See
Figure \ref{fig 3c teo 1}.

\begin{figure}[!h]
\begin{center}\psfrag{A}{\hspace{-.7cm}$L_3<\lambda < L_2$} \psfrag{B}{$\lambda= L_2$} \psfrag{C}{\hspace{.2cm}$L_2 < \lambda
< \beta$} \psfrag{D}{$ \lambda = \beta$} \psfrag{E}{$\lambda >
\beta$}\psfrag{0}{$15_1$} \psfrag{1}{$16_1$} \psfrag{2}{$17_1$}
\psfrag{3}{$18_1$}\psfrag{4}{$19_1$} \epsfxsize=10cm
\epsfbox{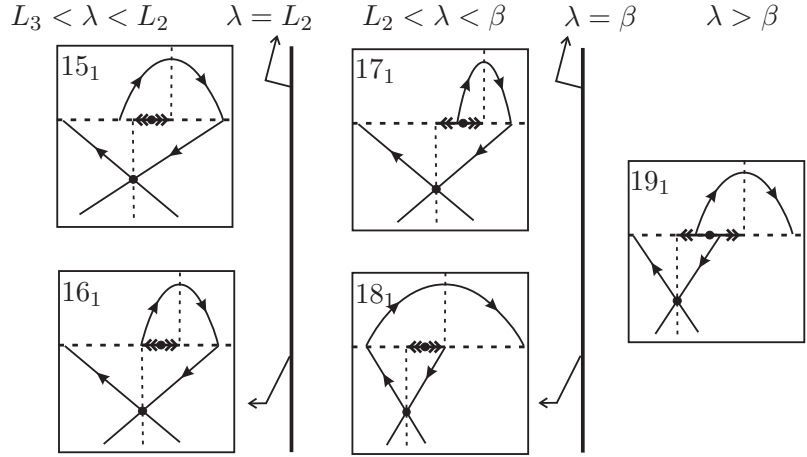} \caption{\small{Cases $15_1 - 19_1$.}}
\label{fig 3c teo 1}
\end{center}
\end{figure}


\begin{figure}[!h]
\begin{center}\psfrag{0}{$\lambda$} \psfrag{1}{$\beta$}\psfrag{A}{$1_1$} \psfrag{B}{$2_1$} \psfrag{C}{$3_1$} \psfrag{D}{$4_1$} \psfrag{E}{$5_1$}
\psfrag{F}{$6_1$}  \psfrag{G}{$7_1$} \psfrag{H}{$8_1$}
\psfrag{I}{$9_1$} \psfrag{J}{$10_1$} \psfrag{K}{$11_1$}
\psfrag{L}{$12_1$}\psfrag{M}{$13_1$} \psfrag{N}{$14_1$}
\psfrag{O}{$15_1$} \psfrag{P}{$16_1$}
\psfrag{Q}{$17_1$}\psfrag{R}{$18_1$} \psfrag{S}{$19_1$}
\epsfxsize=10cm \epsfbox{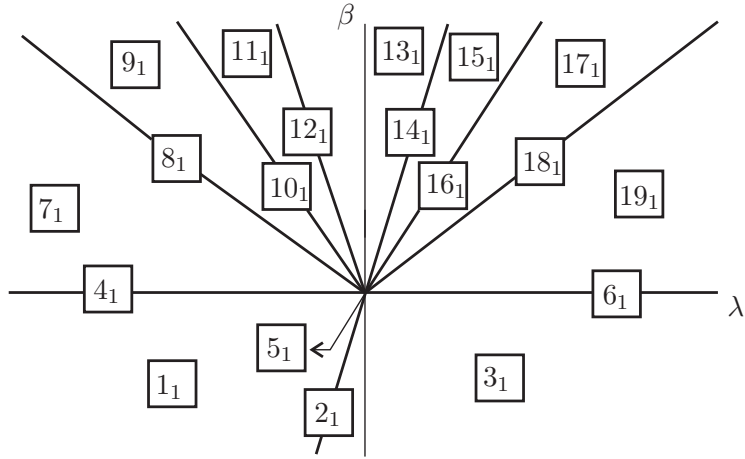}
\caption{\small{Bifurcation Diagram of Theorem 1.}} \label{fig
diagrama bif teo 1}
\end{center}
\end{figure}

The bifurcation diagram is illustrated in Figure \ref{fig diagrama
bif teo 1}.\end{proof}

\begin{remark}\label{observacao separatrizes} In Cases $11_1$ and $15_1$ the ST-bifurcations (as described in \cite{Marcel})
arise. In fact, note that the trajectory passing through $h$ can
make more and more turns around $P$. This fact characterizes a
global bifurcation also reached in other cases. 
\end{remark}


\section{Proof of Theorem 2}\label{secao prova teorema 2}

\begin{proof}[Proof of Theorem 2] In Cases $1_2$, $2_2$ and $3_2$ we assume
that $Y$ presents the behavior $Y^-$.  In Cases $4_2$, $5_2$ and
$6_2$ we assume that $Y$ presents the behavior $Y^0$. In Cases $7_2
- 21_2$ we assume that $Y$ presents the behavior $Y^+$.

$\diamond$ \textit{Cases $(1_2)$ $d_1<e_1$, $(2_2)$ $d_1=e_1$,
$(3_2)$ $d_1>e_1$, $(4_2)$ $d_1<s_1$, $(5_2)$ $d_1=s_1$ and $(6_2)$
$d_1>s_1$:} The analysis of these cases are done in a similar way as
the cases
 $1_1$, $2_1$, $3_1$, $4_1$, $5_1$
and $6_1$.

Observe that as the parameter $\lambda$ increases $-$ assuming the
values $M_0$, $M_1$ and $M_2$ below described $-$ it appear
orbit-arcs of $X$ connecting the points $h$ and $i$, $h$ and $j$ and
$i$ and $j$ respectively. The values of $M_0$, $M_1$ and $M_2$ are:

%

$\begin{array}{ccl}
  M_0 & = & (-3 - 3 \alpha (-2 + \alpha + 2 (-1 +
\alpha) \beta) + \\ && +\sqrt{9 (-1 + \alpha)^4 - 12 (-1 +
\alpha)^2 \beta^2})/(6 (-1 + \alpha)^2), \\
  M_1 & = & -1/2  +\sqrt{9 - 12
\beta^2}/6  \\
  M_2 & = & (-3 + 6 \beta - 3
\alpha (-2 + \alpha + 2 \beta) + \\ && +\sqrt{ 9 (-1 + \alpha)^4 -
12 (-1 + \alpha)^2 \alpha^2 \beta^2})/(6 (-1 +
\alpha)^2).   \\
\end{array}$

$\diamond$ \textit{Cases $(7_2)$ $\lambda< - \beta$, $(8_2)$
$\lambda = - \beta$, $(9_2)$ $- \beta < \lambda < M_0$, $(10_2)$
$\lambda = M_0$ and $(11_2)$ $M_0 < \lambda < M_1 $:} Analogous to
Cases $7_1 - 11_1$ changing $L_0$ by $M_0$ and $L_1$ by
$M_1$. 

$\diamond$ \textit{Case $(12_2)$ $\lambda = M_1$:} The points of
$\Sigma$ outside the interval $(d_1,i_1)$ belong to $\Sigma_c$ and
the points inside this interval belong to $\Sigma_s$. The direction
function $H$ assumes positive values in a neighborhood of $d_1$,
negative values  in a  neighborhood of $i_1$ (see Remark \ref{obs H
positiva ou negativa nas dobras}) and there exists a unique value
$\widetilde{P}=\widetilde{P}_{\lambda,\alpha,\beta}$ such that
$H(\widetilde{P})=0$. So $P=(\widetilde{P},0)$ is a
$\Sigma-$attractor. Since
 $\lambda=M_1$, there is an orbit-arc
$\gamma_{1}^X$ of $X$ connecting the points $h$ and $j$. It
generates a $\Sigma-$graph $\Gamma = \gamma_{1}^{X} \cup \Sigma_e
\cup S \cup \Sigma_c$ of kind I. Since $ \alpha > \alpha_0$, where
$\alpha_0$ is given by  (\ref{eq alpha zero}), it is straightforward
to show that the \textit{First Return Map} 
defined in the interval $(h_1,d_1) \subset \Sigma$ do not have fixed
points. By consequence, $\Gamma$ is a repeller for the trajectories
inside it and canard cycles do not arise. See Figure \ref{fig 1 teo
2}.

\begin{figure}[!h]
\begin{center}\psfrag{A}{$\lambda =M_1$} \psfrag{B}{$M_1<\lambda < i_1$} \psfrag{C}{$\lambda = i_1$} \psfrag{D}{$ \lambda = -\frac{\beta}{2}$}
\psfrag{E}{$-\frac{\beta}{2} < \lambda < 0$}\psfrag{0}{$12_2$}
\psfrag{1}{$13_2$} \psfrag{2}{$14_2$}
\psfrag{3}{$10_1$}\psfrag{4}{$11_1$} \epsfxsize=12.5cm
\epsfbox{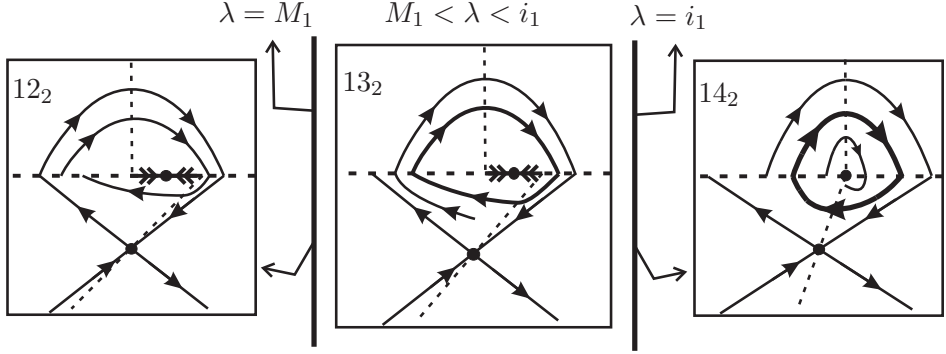} \caption{\small{Cases $12_2$, $13_2$ and
$14_2$.}} \label{fig 1 teo 2}
\end{center}
\end{figure}

$\diamond$ \textit{Case $(13_2)$ $M_1<\lambda < i_1$:} The
distribution of the connected components of $\Sigma$ and the
behavior of $H$ are the same as Case $12_2$. 
Since $M_1<\lambda < i_1$, there is an orbit-arc $\gamma_{1}^X$ of
$X$ connecting $j$ to a point $k=(k_1,0) \in \Sigma$, where $k_1 \in
(h_1,d_1)$, for negative time. Also there is an orbit-arc
$\gamma_{1}^{Y}$ of $Y$ connecting  $k$ to a point $l=(l_1,0) \in
\Sigma$, where $l_1 \in (i_1,j_1)$, for negative time. Repeating
this argument, we can find an increasing sequence $(k_i)_{i \in
\N}$. We can prove that there is an interval $I\subset(k,d)$ such
that $\eta'=(\varphi_{Y}\circ \varphi_{X})' < 1$ on $I$. As $P$ is a
$\Sigma-$attractor, there is an interval $J \subset (k,d)$ such that
$\eta' > 1$ on $J$. Moreover, we can prove that $\eta$ has a unique
fixed point $Q \in (k,d)$. As consequence, there passes through $Q$
a repeller canard cycle $\Gamma$ of kind I. See Figure \ref{fig 1
teo 2}. This canard cycle is born from the bifurcation of the
$\Sigma-$graph present in Case $12_2$. The expression of $\eta$ is
too large, so the general case will be omitted. For the particular
case when $\alpha=-1$, $\beta=1/2$ and $\lambda=-1/2 + 11
\sqrt{6}/60$, the map $\eta$ is given by
$$\begin{array}{ccl}
                                                  \eta(x) & = & \frac{3}{4} +
\frac{3}{2} \Big( - \frac{1}{2} + \frac{11}{10 \sqrt{6}} \Big) +
\frac{x}{2} + \\
                                                   & & -
 \frac{1}{4} \sqrt{3} \sqrt{\Big( 1 - 2 \Big( - \frac{1}{2} + \frac{11}{10 \sqrt{6}} \Big) - 2 x \Big) \Big(3 +
    2 \Big(- \frac{1}{2} + \frac{11}{10 \sqrt{6}} \Big) + 2 x \Big)}
                                                 \end{array}.$$
A straightforward calculation shows that the unique fixed point of
this par\-ti\-cu\-lar $\eta$ occurs when $x = -\sqrt{29/2}/10$.

$\diamond$ \textit{Case $(14_2)$ $\lambda = i_1$:} Every point of
$\Sigma$ belongs to $\Sigma_c$ except the point $d=i$. The canard
cycle presented in the previous case is persistent for this case
(remember that this canard cycle is born from the bifurcation of the
$\Sigma-$graph of Case $12_2$. So, it radius does not tend to zero
when $\lambda$ tends to $i_1$). So the non generic tangential
singularity $d=i$ behaves itself like a weak attractor focus. See
Figure \ref{fig 1 teo 2}.

$\diamond$ \textit{Cases $(15_2)$ $i_1< \lambda < M_3$ and $(16_2)$
$\lambda = M_3$:} Analogous to Cases $13_1 - 14_1$ replacing $L_1$
by $i_1$ and $L_3$ by $M_3$, where $M_3$ is the limit value for
which $\Gamma_1$ collides with $\Gamma_2$.

$\diamond$ \textit{Cases $(17_2)$ $M_3<\lambda< M_2$, $(18_2)$
$\lambda = M_2$, $(19_2)$ $M_2 < \lambda < \beta$, $(20_2)$ $\lambda
= \beta$ and $(21_2)$ $\lambda
> \beta$:}  Analogous to
Cases $15_1 - 19_1$ replacing $L_2$ by $M_2$ and $L_3$ by $M_3$. 

\begin{figure}[!h]
\begin{center}\psfrag{0}{$\lambda$} \psfrag{1}{$\beta$}
\psfrag{A}{$1$} \psfrag{B}{$2$} \psfrag{C}{$3$} \psfrag{D}{$4$}
\psfrag{E}{$5$} \psfrag{F}{$6$}  \psfrag{G}{$7$} \psfrag{H}{$10$}
\psfrag{I}{$11$} \psfrag{J}{$12$} \psfrag{K}{$13$}
\psfrag{L}{$14$}\psfrag{M}{$15$} \psfrag{N}{$16$} \psfrag{O}{$17$}
\psfrag{P}{$18$} \psfrag{Q}{$19$}\psfrag{R}{$20$}
\psfrag{S}{$21$}\psfrag{2}{$8$} \psfrag{3}{$9$} \epsfxsize=10cm
\epsfbox{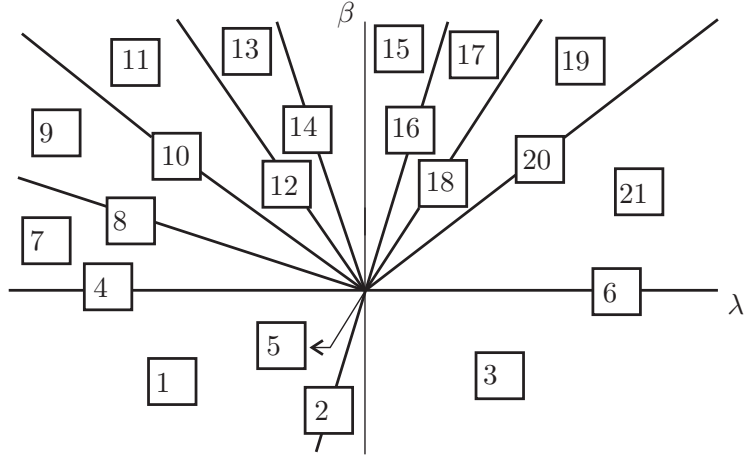} \caption{\small{Bifurcation Diagram
of Theorems 2 and 3.}} \label{fig diagrama bif teo 2}
\end{center}
\end{figure}

The bifurcation diagram is illustrated in Figure \ref{fig diagrama
bif teo 2}.\end{proof}


\section{Proof of Theorem 3}\label{secao prova teorema 3}

\begin{proof}[Proof of Theorem 3] In Cases $1_3$, $2_3$ and $3_3$ we assume
that $Y$ presents the behavior $Y^-$.  In Cases $4_3$, $5_3$ and
$6_3$ we assume that $Y$ presents the behavior $Y^0$. In Cases $7_3
- 21_3$ we assume that $Y$ presents the behavior $Y^+$.

$\diamond$ \textit{Cases $(1_3)$ $d_1<e_1$, $(2_3)$ $d_1=e_1$,
$(3_3)$ $d_1>e_1$, $(4_3)$ $d_1<s_1$, $(5_3)$ $d_1=s_1$ and $(6_3)$
$d_1>s_1$:} Analogous to Cases $1_1$, $2_1$, $3_1$, $4_1$, $5_1$ and
$6_1$.

In what follows we consider $M_0$, $M_1$, $M_2$ and $M_3$ as in the
previous theorem.

$\diamond$ \textit{Cases $(7_3)$ $\lambda< - \beta$, $(8_3)$
$\lambda = - \beta$, $(9_3)$ $- \beta < \lambda < M_0$, $(10_3)$
$\lambda = M_0$ and $(11_3)$ $M_0 < \lambda < i_1 $:}
Analogous to Cases $7_2 - 11_2$ changing $M_1$ by $i_1$.

$\diamond$ \textit{Case $(12_3)$ $\lambda = i_1$:} Every point of
$\Sigma/\{d\}$ belongs to $\Sigma_c$. In a similar way as Case
$13_2$, we can construct sequences $(k_i)_{i \in \N}$ and $(l_i)_{i
\in \N}$. Since $d=i$ we have that $k_i \rightarrow d$ and $l_i
\rightarrow d$. So $d$ is a non generic tangential singularity that
behaves itself like an attractor. See Figure \ref{fig 1 teo 3}.

\begin{figure}[!h]
\begin{center}\psfrag{A}{$\lambda = i_1$} \psfrag{B}{$i_1<\lambda < M_1$} \psfrag{C}{$\lambda =M_1$}
\psfrag{1}{$13_3$} \psfrag{2}{$14_3$} \psfrag{0}{$12_3$}
\epsfxsize=12.5cm \epsfbox{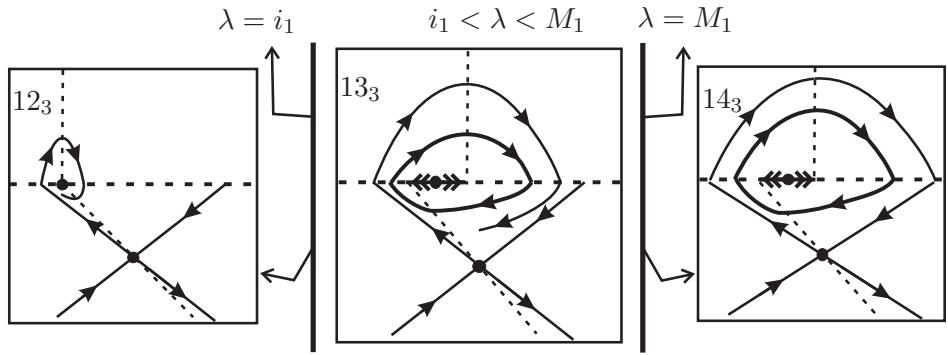} \caption{\small{Cases
$12_3$, $13_3$ and $14_3$.}} \label{fig 1 teo 3}
\end{center}
\end{figure}

$\diamond$ \textit{Case $(13_3)$ $i_1 <\lambda < M_1$:} Analogous to
Case $13_2$ except that there is a change of stability on
$P=(\widetilde{P},0)$
, which is a $\Sigma-$repeller, and on $\Gamma$, which is an
attractor canard cycle of kind I. This canard cycle is born from the
bifurcation of the non-generic tangential singularity of Case
$12_3$. See Figure \ref{fig 1 teo 3}.

$\diamond$ \textit{Case $(14_3)$ $\lambda = M_1$:} Analogous to Case
$12_2$ except that occurs a change of stability on $P=(\widetilde{P},0)$, which is a $\Sigma-$repeller. 
 This fact generates a bifurcation like Hopf near $P$ and there appears a hyperbolic attractor canard cycle $\Gamma_1$, of kind I, between $P$ and the
 $\Sigma-$graph $\Gamma_2$. See Figure \ref{fig 1 teo
3}.

$\diamond$ \textit{Cases $(15_3)$ $M_1<\lambda< M_3$ and $(16_3)$
$\lambda = M_3$:}  Analogous to Cases $15_2 - 16_2$, replacing $i_1$
by $M_1$.

$\diamond$ \textit{Cases $(17_3)$ $M_3<\lambda< M_2$, $(18_3)$
$\lambda = M_2$, $(19_3)$ $M_2 < \lambda < \beta$, $(20_3)$ $\lambda
= \beta$ and $(21_2)$ $\lambda
> \beta$:}  Analogous to
Cases $17_2 - 21_2$. 


The bifurcation diagram is illustrated in Figure \ref{fig diagrama
bif teo 2}.\end{proof}

\section{Proof of Theorem 4}\label{secao prova teorema 4}


\begin{proof}[Proof of Theorem 4] Since $X$ has a unique $\Sigma-$fold point which is visible we conclude that canard cycles do not arise.

In Cases $1_4$, $2_4$ and $3_4$ we assume that $Y$ presents the
behavior $Y^-$. In Cases $4_4$, $5_4$ and $6_4$ we assume that $Y$
presents the behavior $Y^0$. In these cases, when
 it is well defined, the direction function $H$ assumes positive values.

$\diamond$ \textit{Case $(1_4)$ $d_1<e_1$:} The points of $\Sigma$
inside the interval $(d_1,e_1)$ belong to $\Sigma_c$. The points on
the left of
 $d_1$ belong to $\Sigma_s$ and the points on the right of
 $e_1$ belong to $\Sigma_e$. See Figure \ref{fig 1 teo 4}.

$\diamond$ \textit{Case $(2_4)$ $d_1=e_1$:} Here $\Sigma_c =
\emptyset$.
  The vector fields $X$ and $Y$ are linearly
 dependent on $d_1=e_1$ which is a tangential singularity. Moreover, it is an attractor
 for the trajectories of $Z$ crossing $\Sigma_s$ and  a repeller for the trajectories of $Z$ crossing
 $\Sigma_e$. See Figure \ref{fig 1 teo 4}.

$\diamond$ \textit{Case $(3_4)$ $d_1>e_1$:} The points of $\Sigma$
inside the interval $(e_1,d_1)$ belong to $\Sigma_c$. The points on
the left of
 $e_1$ belong to $\Sigma_s$ and the points on the right of
 $d_1$ belong to $\Sigma_e$. See Figure \ref{fig 1 teo 4}.

\begin{figure}[!h]
\begin{center}\psfrag{A}{$\lambda < 0$} \psfrag{B}{$\lambda=0$} \psfrag{C}{$\lambda > 0$}\psfrag{0}{$1_4$} \psfrag{1}{$2_4$} \psfrag{2}{$3_4$}
\epsfxsize=9.5cm \epsfbox{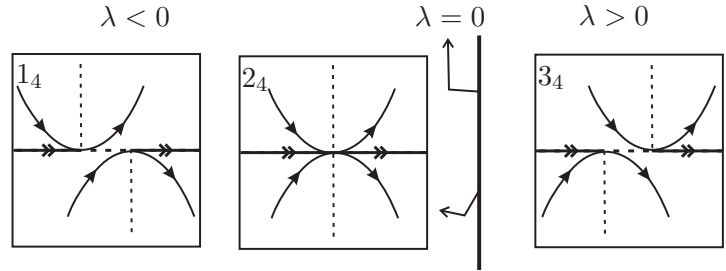} \caption{\small{Cases
$1_4$, $2_4$ and $3_4$.}} \label{fig 1 teo 4}
\end{center}
\end{figure}


$\diamond$ \textit{Case $(4_4)$ $d_1<s_1$:} The points of $\Sigma$
inside the interval $(d_1,s_1)$ belong to $\Sigma_c$. The points on
the left of
 $d_1$ belong to $\Sigma_s$ and the points on the right of
 $s_1$ belong to $\Sigma_e$. See Figure \ref{fig 2 teo 4}.

$\diamond$ \textit{Case $(5_4)$ $d_1=s_1$:} Here
$\Sigma_c=\emptyset$ and $S$ is an attractor
 for the trajectories of $Z$ crossing $\Sigma_s$ and it is a repeller for the trajectories of $Z$ crossing
 $\Sigma_e$. See Figure \ref{fig 2 teo 4}.

$\diamond$ \textit{Case $(6_4)$ $d_1>s_1$:} The points of $\Sigma$
inside the interval $(d_1,s_1)$ belong to $\Sigma_c$. The points on
the left of
 $s_1$ belong to $\Sigma_s$ and the points on the right of
 $d_1$ belong to $\Sigma_e$. See Figure \ref{fig 2 teo 4}.

\begin{figure}[!h]
\begin{center}\psfrag{A}{\hspace{-.5cm}$\lambda < 0$} \psfrag{B}{$\lambda=0$} \psfrag{C}{$\lambda > 0$}\psfrag{0}{$4_4$} \psfrag{1}{$5_4$} \psfrag{2}{$6_4$}
\epsfxsize=9.5cm \epsfbox{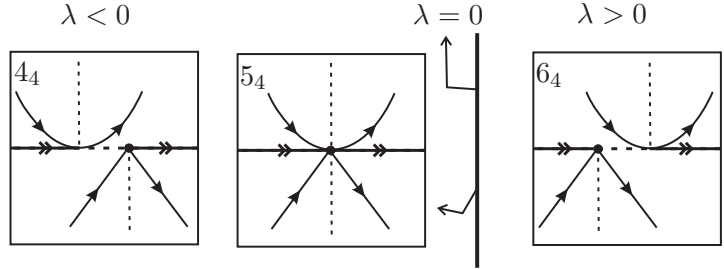} \caption{\small{Cases
$4_4$, $5_4$ and $6_4$.}} \label{fig 2 teo 4}
\end{center}
\end{figure}

In Cases $7_4 - 13_4$ we assume that $Y$ presents the behavior
$Y^+$.

$\diamond$ \textit{Cases $(7_4)$ $d_1 < h_1$, $(8_4)$ $d_1 = h_1$
and $(9_4)$ $h_1 < d_1 < i_1$:} The points of $\Sigma$ inside the
interval $(d_1,i_1)$ belong to $\Sigma_c$. The points on the left of
 $d_1$ belong to $\Sigma_s$ and the points on the right of
 $i_1$ belong to $\Sigma_e$. The direction function $H$ assumes positive
values on $\Sigma_s$ and negative values in a neighborhood of $i_1$.
Moreover, $H(\beta \lambda / (-1 + \beta))=0$ and the
$\Sigma-$repeller $P=(\beta \lambda / (-1 + \beta),0)$ is the unique
pseudo equilibrium. See Figure \ref{fig 3a teo 4}.

$\diamond$ \textit{Case $(10_4)$ $d_1=i_1$:} Here $\Sigma_c =
\emptyset$. The vector fields $X$ and $Y$ are linearly
 dependent on the tangential singularity $d_1=i_1$.
A straightforward calculation shows that $H(z)=(1-\beta)/2 \neq 0$
for all $z \in \Sigma/\{ d \}$. So $d_1=i_1$ is an attractor
 for the trajectories of $Z$ crossing $\Sigma_s$ and  a repeller for the trajectories of $Z$ crossing
 $\Sigma_e$. Moreover, $\Delta=\{ d \} \cup \overline{d j} \cup \Sigma_e \cup \{ S \} \cup \Sigma_c \cup \overline{h d}$
 is a $\Sigma-$graph of kind III in such a way that each $Q$ in its interior belongs to another $\Sigma-$graph of kind III passing through $d$.  See Figure
\ref{fig 3a teo 4}.

\begin{figure}[!h]
\begin{center}\psfrag{A}{$\lambda < -\beta$} \psfrag{B}{$\lambda = - \beta$} \psfrag{C}{$- \beta < \lambda
< 0$} \psfrag{D}{$ \lambda = 0$}\psfrag{0}{$7_4$} \psfrag{1}{$8_4$}
\psfrag{2}{$9_4$} \psfrag{3}{$10_4$} \epsfxsize=12.5cm
\epsfbox{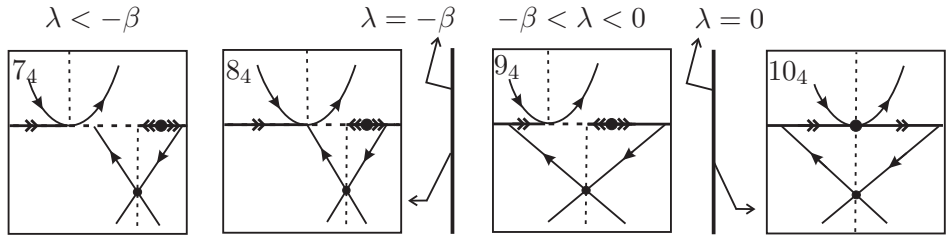} \caption{\small{Cases $7_4 - 10_4$.}}
\label{fig 3a teo 4}
\end{center}
\end{figure}

$\diamond$ \textit{Cases $(11_4)$ $i_1 < d_1 < j_1$, $(12_4)$ $d_1 =
j_1$ and $(13_4)$ $j_1 < d_1$:} The points of $\Sigma$ inside the
interval $(i_1,d_1)$ belong to $\Sigma_c$. The points on the left of
 $i_1$ belong to $\Sigma_s$ and the points on the right of
 $d_1$ belong to $\Sigma_e$. The direction function $H$ assumes
 positive
values on $\Sigma_e$ and negative values in a neighborhood of $i_1$.
Moreover, $H(\beta \lambda / (-1 + \beta))=0$ and the
$\Sigma-$attractor $P=(\beta \lambda / (-1 + \beta),0)$ is the
unique pseudo equilibrium. See Figure \ref{fig 3b teo 4}.

\begin{figure}[!h]
\begin{center}\psfrag{A}{$0<\lambda < \beta$} \psfrag{B}{$\lambda = \beta$} \psfrag{C}{$\beta < \lambda$}
\psfrag{0}{$11_4$} \psfrag{1}{$12_4$} \psfrag{2}{$13_4$}
 \epsfxsize=9.5cm
\epsfbox{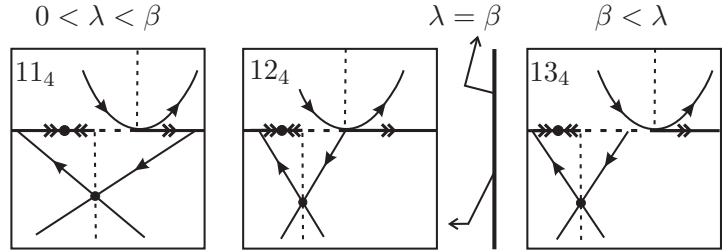} \caption{\small{Cases $11_4 - 13_4$.}}
\label{fig 3b teo 4}
\end{center}
\end{figure}

\begin{figure}[!h]
\begin{center}\psfrag{0}{$\lambda$} \psfrag{1}{$\beta$}\psfrag{A}{\hspace{.15cm}$1$} \psfrag{B}{\hspace{.1cm}$2$}
 \psfrag{C}{\hspace{.15cm}$3$} \psfrag{D}{\hspace{.12cm}$4$} \psfrag{E}{\hspace{.13cm}$5$}
\psfrag{F}{\hspace{.1cm}$6$}  \psfrag{G}{\hspace{.15cm}$7$}
\psfrag{H}{\hspace{.15cm}$8$} \psfrag{I}{\hspace{.17cm}$9$}
\psfrag{J}{\hspace{.15cm}$10$} \psfrag{K}{\hspace{.15cm}$11$}
\psfrag{L}{\hspace{.1cm}$12$}\psfrag{M}{\hspace{.1cm}$13$}
\epsfxsize=7cm \epsfbox{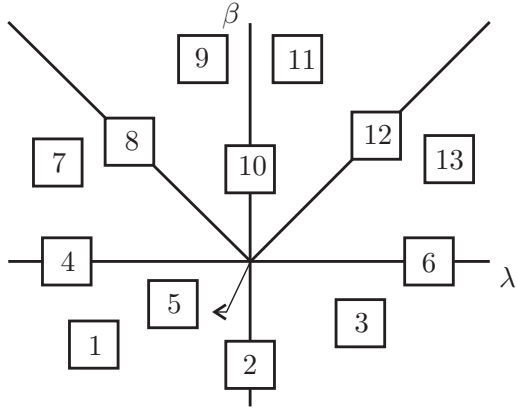}
\caption{\small{Bifurcation Diagram of Theorems 4, 5 and 6.}}
\label{fig diagrama bif teo 4 5 e 6}
\end{center}
\end{figure}

The bifurcation diagram is illustrated in Figure \ref{fig diagrama
bif teo 4 5 e 6}.\end{proof} 


\section{Proof of Theorem 5}\label{secao prova teorema 5}

\begin{proof}[Proof of Theorem 5] The
direction function $H$ has a  root $Q=(q,0)$
where\begin{equation}\label{eq p novo}\begin{array}{cl}q=&
\dfrac{1}{2(\alpha+1)} ((-1+\alpha)(1-\beta)- \lambda(1+\alpha) +\\&
+ \sqrt{((-1+\alpha)(1-\beta)- \lambda(1+\alpha))^2 + 4 \beta
(1+\alpha) (1+ \alpha + \lambda(-1+\alpha))
}).\end{array}\end{equation}Moreover, $H$ assumes positive values on
the right of $Q$ and negative values on the left of $Q$. Note that
when $\alpha \rightarrow -1$ so $Q \rightarrow - \infty$ under the
line $\{y=0 \}$ and it occurs the configurations showed in Theorem
4.

In Cases $1_5$, $2_5$ and $3_5$ we assume that $Y$ presents the
behavior $Y^-$.  In Cases $4_5$, $5_5$ and $6_5$ we assume that $Y$
presents the behavior $Y^0$. In Cases $7_5 - 13_5$ we assume that
$Y$ presents the behavior $Y^+$.

$\diamond$ \textit{Cases $(1_5)$ $d_1<e_1$, $(2_5)$ $d_1=e_1$,
$(3_5)$ $d_1>e_1$, $(4_5)$ $d_1<s_1$, $(5_5)$. $d_1=s_1$ and $(6_5)$
$d_1>s_1$:} Analogous to Cases $1_4$, $2_4$, $3_4$, $4_4$, $5_4$ and
$6_4$ respectively, except that here it appears the $\Sigma-$saddle
$Q$ on the left of $d$ and $e$ or $S$. See Figure \ref{fig 1 teo 5}.

\begin{figure}[!h]
\begin{center}\psfrag{A}{\hspace{-2.6cm}$\lambda < (1+\alpha) \beta/(1-\alpha)$}
\psfrag{B}{\hspace{-2.2cm}$\lambda=(1+\alpha) \beta/(1-\alpha)$}
\psfrag{C}{$\lambda > (1+\alpha) \beta/(1-\alpha)$}
\psfrag{0}{$1_5$} \psfrag{1}{$2_5$} \psfrag{2}{$3_5$}
\epsfxsize=10cm \epsfbox{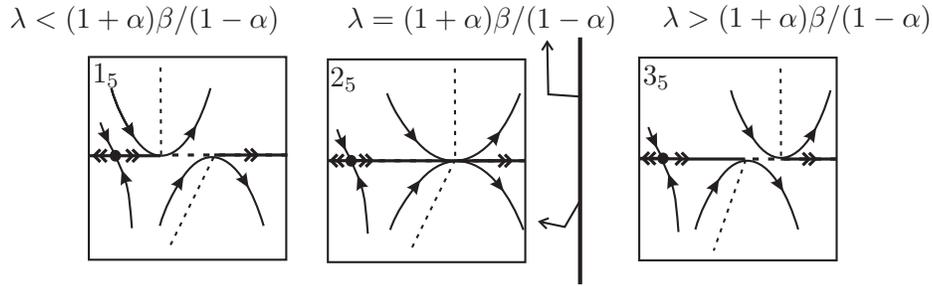} \caption{\small{Cases
$1_5$, $2_5$ and $3_5$.}} \label{fig 1 teo 5}
\end{center}
\end{figure}

%


$\diamond$ \textit{Cases $(7_5)$ $d_1< h_1$, $(8_5)$ $d_1 = h_1$,
$(9_5)$ $h_1 < d_1 < i_1$:} Analogous to Cases $7_4 - 9_4$, except
that here the $\Sigma-$saddle $Q$ appears on the left of $d_1$ and
$i_1$. So $P=(p,0)$ where\begin{equation}\label{eq p novo
2}\begin{array}{cl}p=& \dfrac{1}{2(\alpha+1)} ((-1+\alpha)(1-\beta)-
\lambda(1+\alpha) +\\& - \sqrt{((-1+\alpha)(1-\beta)-
\lambda(1+\alpha))^2 + 4 \beta (1+\alpha) (1+ \alpha +
\lambda(-1+\alpha)) }).\end{array}\end{equation}

$\diamond$ \textit{Case $(10_5)$ $d_1=i_1$:} Analogous to Case
$10_4$, except that here appear the $\Sigma-$saddle $Q$ on the left
of
$d_1=i_1$.

$\diamond$ \textit{Cases $(11_5)$ $i_1 < d_1 < j_1$, $(12_5)$ $d_1 =
j_1$ and $(13_5)$ $j_1 < d_1$:}  Analogous to Cases $11_4 - 13_4$,
except that here the $\Sigma-$saddle $Q$ appears on the left of
$d_1$ and $i_1$.

%

The bifurcation diagram is illustrated in Figure \ref{fig diagrama
bif teo 4 5 e 6}.\end{proof} 


\section{Proof of Theorem 6}\label{secao prova teorema 6}

\begin{proof}[Proof of Theorem 6] The
direction function $H$ has a  root $Q=(q,0)$ where $q$ is given by
 (\ref{eq p novo}).
Moreover, $H$ assumes positive values on the left of $Q$ and
negative values on the right of $Q$. Note that when $\alpha
\rightarrow -1$ so $Q \rightarrow \infty$ under the line $\{y=0 \}$
and  the configurations shown in Theorem 4 occur.

In Cases $1_6$, $2_6$ and $3_6$ we assume that $Y$ presents the
behavior $Y^-$.  In Cases $4_6$, $5_6$ and $6_6$ we assume that $Y$
presents the behavior $Y^0$. In Cases $7_6 - 13_6$ we assume that
$Y$ presents the behavior $Y^+$.

$\diamond$ \textit{Cases $(1_6)$ $d_1<e_1$, $(2_6)$ $d_1=e_1$,
$(3_6)$ $d_1>e_1$, $(4_6)$ $d_1<s_1$, $(5_6)$ $d_1=s_1$ and $(6_6)$
$d_1>s_1$, $(7_6)$ $d_1 < h_1$, $(8_6)$ $d_1 = h_1$, $(9_6)$ $h_1 <
d_1 < i_1$, $(10_6)$ $d_1=i_1$, $(11_6)$ $i_1 < d_1 < j_1$, $(12_6)$
$d_1 = j_1$ and $(13_6)$ $j_1 < d_1$:} Analogous to Cases $1_5$,
$2_5$, $3_5$, $4_5$, $5_5$, $6_5$, $7_5$, $8_5$, $9_5$, $10_5$,
$11_5$, $12_5$ and $13_5$ respectively, except that here the
$\Sigma-$saddle $Q$ takes place on the right of $d_1$, $e_1$, $s_1$
and $i_1$ when these points
appear. 

%

The bifurcation diagram is illustrated in Figure \ref{fig diagrama
bif teo 4 5 e 6}.\end{proof} 

\noindent {\textbf{Acknowledgments.}  The first and the third
authors are partially supported by a FAPESP-BRAZIL grant
2007/06896-5. The second author was partially supported by
FAPESP-BRAZIL grants 2007/08707-5, 2010/18190-2 and 2012/00481-6.}


\end{document}